\documentclass[11pt]{amsart}
\allowdisplaybreaks[1]

\usepackage{a4wide}
\usepackage{amssymb}
\usepackage{amsmath}
\usepackage{amscd} 

\newtheorem{theorem}{Theorem}
\newtheorem{lemma}{Lemma}
\newtheorem{prop}{Proposition}
\newtheorem{coro}{Corollary}
\newtheorem{fact}{Fact}
\newtheorem{definition}{Definition}

\newcommand{\ts}{\hspace{0.5pt}}

\newcommand{\CC}{\mathbb{C}\ts}
\newcommand{\RR}{\mathbb{R}\ts}
\newcommand{\QQ}{\mathbb{Q}\ts}
\newcommand{\ZZ}{\mathbb{Z}}
\newcommand{\NN}{\mathbb{N}}

\newcommand{\oplam}{\mbox{\Large $\curlywedge$}}

\newcommand{\dd}{\,{\rm d}}

\newcommand{\UKV}{U_{K,V}}

\newcommand{\MM}{\mathcal{M}(G)}
\newcommand{\MCV}{\mathcal{M}_{C,V}(G)}
\newcommand{\MTB}{\mathcal{M}^{\infty}(G)}

\newcommand{\cMO}{\mathcal{M}(\varOmega)}
\newcommand{\cMT}{\mathcal{M}(\varTheta)}
\newcommand{\cPO}{\mathcal{P}(\varOmega)}
\newcommand{\cPGO}{\mathcal{P}_G (\varOmega)}
\newcommand{\cPGT}{\mathcal{P}_G (\varTheta)}

\newcommand{\Oomega}{(\varOmega,\alpha)}
\newcommand{\Oop}{(\varOmega_{\rm p},\alpha)}
\newcommand{\Ttheta}{(\varTheta,\beta)}

\newcommand{\LO}{L^2 (\varOmega,m)}
\newcommand{\LT}{L^2 (\varTheta,n)}

\newcommand{\cV}{\mathcal{V}}

\newcommand{\cDV}{\mathcal{D}_V (G)}

\newcommand{\cUD}{\ts\mathcal{UD}(G)}

\newcommand{\Ghat}{\widehat{G}}
\newcommand{\gammahat}{\widehat{\gamma}}

\newcommand{\supp}{\mbox{supp}}

\begin{document}

\title[Delone dynamical systems]
{Deformation of Delone dynamical systems \\[1mm] 
and pure point diffraction}

\author{Michael Baake}
\address{Fakult\"at f\"{u}r Mathematik, Universit\"{a}t Bielefeld,
Postfach 100131, 33501 Bielefeld, Germany}
\email{mbaake@math.uni-bielefeld.de}

\author{Daniel Lenz}
\address{Fakult\"at f\"ur Mathematik, TU Chemnitz,
09107 Chemnitz, Germany}
\email{dlenz@mathematik.tu-chemnitz.de}

\begin{abstract} 
This paper deals with certain dynamical systems built from point sets
and, more generally, measures on locally compact Abelian
groups. These systems arise in the study of quasicrystals and
aperiodic order, and important subclasses of them exhibit pure point 
diffraction spectra. We discuss the relevant framework and  recall fundamental
results and examples. In particular, we show that pure point diffraction is 
stable under ``equivariant'' local perturbations and discuss various examples,
including deformed model sets. 
A key step in the proof of stability consists in transforming the
problem into a question on factors of dynamical systems.
\end{abstract}

\maketitle

\section{Introduction}

Aperiodic order has become a topic of intense research over the last two
decades \cite{Mbook,Pat,BMbook,Suck,Trebin}. While the term is not rigorously 
defined (yet), it roughly refers to forms of order at the very verge between 
periodic and non-periodic structures. As such, it has attracted attention in 
various branches of mathematics including geometry, combinatorics, ergodic theory, 
operator theory and harmonic analysis.
  
An important trigger in these developments has been the actual discovery of
physical substances with strong aperiodic order \cite{dany}, which are 
now called quasicrystals. They owe their discovery to their remarkable diffraction
patterns: These patterns imply a high degree of order as they are pure point
spectra (or Bragg spectra), while, at the same time, they exclude periodicity 
by their non-crystallographic symmetries. Accordingly, the study of pure point 
diffraction has been an important topic in this context ever since.
  
This paper is concerned with pure point diffraction. More precisely, we
study the stability of pure point diffraction under certain deformations.
This issue is a very natural one, both from the physical and the
mathematical point of view. In order to study stability under deformations, 
we need to review the undeformed case first. To make the paper essentially 
self-contained, this discussion is carried out at some length, including 
relevant concepts and examples. Moreover, we hope that 
the paper can serve as an introductory survey on the treatment of
diffraction via dynamical systems for the reader unfamiliar with the field.

\smallskip

Delone sets provide an important model class for the description of
aperiodic order. In particular, they can be viewed as a mathematical
abstraction of the set of atomic positions of a physical quasicrystal
(at zero temperature, or at a given instant of time). Many of the
rather intriguing spectral properties of quasicrystals can be
formulated, in a simplified manner, on the basis of Delone sets.
This is also a rather common class of structures in the mathematical
theory of aperiodic order \cite{Lag}. It is attractive because it admits a direct
geometric interpretation with two Delone sets being close to one another
if large patches (around some fixed point of reference, say) coincide,
possibly after a tiny local rearrangement of the individual points.

However, from a more physical point of view, other scenarios are also
very important. In particular, the description of an aperiodic distribution
of matter by means of (continuous) quasi- or almost periodic density
functions has been emphasized right from the very beginning of quasicrystal
theory \cite{Bak}. Here, closeness of two structures is more adequately
described by means of the supremum norm, as in the theory of almost
periodic functions.

As is apparent, these two pictures are not compatible --- unless they are
embedded into a larger class of structures that admit both the Delone
(or tiling) picture and the continuous density description as special cases.
One possibility is the use of translation bounded (complex)
measures, equipped with the vague topology. Here, two structures (i.e.,
measures) are close if their evaluations with continuous functions
supported on a large compact set $K$ are close. This entails both
situations mentioned above, one being described by pure point measures,
the other by absolutely continuous measures with continuous
Radon-Nikodym densities.

In view of the fact that the original distinction between the discrete
and the continuous approach led to rather hefty disputes on the
justification and appropriateness of the two approaches, we believe that
the systematic development of a unified frame is overdue. We take
this as our main motivation for a dynamical systems approach based on
{\em measures}, though we shall also spell out the details for the more
conventional (and perhaps more intuitive) approach via Delone sets.

\smallskip

As mentioned already, one important issue in this context is that of the
{\em stability\/} of certain features, e.g., stability under slight modifications 
or deformations.  The question of stability of pure point diffractivity is
addressed in this paper.

Our main abstract result shows that pure point diffraction is stable under local
``equivariant'' perturbations. The proof relies on two steps: We use a
recent result of ours \cite{BL} (see \cite{LMS-1,Gouere-2} for related
material) which establishes when pure point dynamical spectrum is equivalent 
to pure point diffraction spectrum. This effectively transforms the
stability problem into a question on dynamical systems. This question
is then solved by studying certain factors of the original dynamical
system.

To give the reader a flavour of this procedure, we include the following 
rather informal statement of our main result, when restricted to Delone sets.

\smallskip \noindent
{\bf Result.} {\it The hull of an admissibly deformed Delone set is a topological 
factor of the hull of the original Delone set. In particular, if a Delone set has 
pure point diffraction spectrum, its deformation has pure point 
diffraction spectrum as well. }

\smallskip

A precise version of this result is given in Theorem \ref{stabilitypoint}. As 
mentioned already, our setting is general enough to treat not only the case of 
Delone sets but rather the case of arbitrary measure dynamical systems. This is 
made precise in Theorem \ref{stabilityabstract}. 

\smallskip

The abstract result is applied to various examples. In particular, we
study perturbations of model sets in the context of cut and project schemes. 
This generalizes the corresponding considerations of Hof \cite{Hof3} and Bernuau 
and Duneau \cite{BD}. It also shows that related results of Clark and Sadun
\cite{CS} fall well within our framework.

\smallskip

Our results should be compared to complementary results of Hof \cite{Hof4}. 
They show that {\em random}\/ perturbations do not leave a pure point
spectrum unchanged, but rather introduce an absolutely continuous
component, see also \cite{BH,BM1} for further examples.

We are well aware of the fact that considerable parts of the following
investigation dealing with topological dynamical systems can be
generalized to measurable dynamical systems.  However, by its very
nature, the subject of aperiodic order seems to be a topological
one. For this reason, we stick to the topological category.

\smallskip

The paper is organized as follows. In Section \ref{Generalities}, we
introduce some basic notation concerning topological dynamical systems. In
Section \ref{factor}, we recall and establish various facts on factors. These
considerations are the abstract core behind our deformation procedure. 
Section \ref{Measure} is devoted to a discussion of diffraction in the context  
of dynamical systems of Delone sets and measures. 
The abstract deformation procedure and the stability of pure
point diffraction under this type of deformation is discussed in
Section \ref{Framework}. Applications to model sets are studied in
Section \ref{Modelset}, which also contains a brief summary of their
general definition. The various concepts and results will then be
illustrated with a concrete example, the silver mean chain, in Section
\ref{example}. Further aspects of the deformation procedure, in particular
concerning topological conjugacy, are discussed in Section \ref{Further}.

\section{Generalities on dynamical systems} 
\label{Generalities}
Our considerations are set in the framework of topological dynamical
systems. We are dealing with $\sigma$-compact locally compact topological 
groups and compact spaces. Thus, we start with some basic notation 
and facts concerning locally compact topological spaces used 
throughout the paper.

Whenever $X$ is a $\sigma$-compact locally compact space 
(by which we mean to include
the Hausdorff property), we denote the space of continuous functions on
$X$ by $C(X)$ and the subspace of continuous functions with compact
support by $C_c (X)$. This space is equipped with the locally convex
limit topology induced by the canonical embeddings
$C_K(X)\hookrightarrow C_c (X)$, where $C_K (X)$ is the space of
complex continuous functions with support in a given compact set 
$K\subset X$. Here, each $C_K (X)$ is equipped with the topology
induced by the standard supremum norm.

As $X$ is a topological space, it carries a natural $\sigma$-algebra,
namely the Borel $\sigma$-algebra generated by all closed subsets of
$X$.  The set $\mathcal{M} (X)$ of all complex regular Borel measures
on $G$ can then be identified with the space $C_c (X)^\ast$ of
complex valued, continuous linear functionals on $C_c(G)$. This is justified
by the Riesz-Markov representation theorem, compare \cite[Ch.\ 6.5]{Ped}
for details. In particular, we can
write $\int_X f \dd\mu = \mu(f)$ for $f\in C_c(X)$ and simplify the 
notation this way. The space $\mathcal{M} (X)$
carries the vague topology, i.e., the weakest topology that makes all
functionals $\mu\mapsto \mu(\varphi)$, $\varphi\in C_c (X)$,
continuous.  The total variation of a measure $\mu \in \mathcal{M} (X)$ 
is denoted by $|\mu|$.

\smallskip

We now fix a $\sigma$-compact locally compact Abelian (LCA) group $G$ 
for the remainder of the paper. The dual group of $G$ is denoted by
$\Ghat$, and the pairing between a character $\hat{s} \in \Ghat$ and $t
\in G$ is written as $(\hat{s},t)$.  Whenever $G$ acts on the compact
space $\varOmega$ (which is then also Hausdorff by our
convention) by a continuous action
\begin{equation*}
   \alpha \! : \; G\times \varOmega \; \longrightarrow \; \varOmega
   \, , \quad (t,\omega) \, \mapsto \, \alpha^{}_{t} (\omega) \, ,
\end{equation*} 
where $G\times \varOmega$ carries the product topology, the pair 
$\Oomega$ is called a {\em topological dynamical system\/} over $G$. We
shall often write $\alpha^{}_{t}\ts \omega$ for $\alpha^{}_{t} (\omega)$.
If $\omega\in\varOmega$ satisfies $\alpha^{}_{t}\ts\omega = \omega$, $t$
is called a {\em period\/} of $\omega$. If all $t\in G$ are periods,
$\omega$ is called $G$-invariant, or $\alpha$-{\em invariant\/} to 
refer to the action involved.

The set of all Borel probability measures on $\varOmega$ is denoted by 
$\cPO$, and the subset of $\alpha$-invariant probability measures by $\cPGO$.  
As $\varOmega$ is compact,  $C_c(\varOmega)$ equipped with the supremum 
norm is a Banach space. The vague topology on $\mathcal{M} (\varOmega)$ is 
then just the weak-$\ast$ topology. By Alaoglu's theorem on weak-$\ast$
compactness of the unit sphere (compare \cite[Thm.\ 2.5.2]{Ped}), we 
easily conclude that $\cPO$ is compact.  As
$\cPGO$ is obviously closed in $\cPO$, it is then  compact as well. 
Apparently, $\cPGO$ is convex. More importantly, it is always
non-empty. For $G=\ZZ$, this is standard, compare  Section 6.2 in \cite{Wal}.  
This proof only uses the existence of a van Hove sequence
and the compactness of $\cPO$. Thus, it can be carried over to
our setting (for the existence of van Hove sequences, we refer the reader to
\cite[p.~145]{Martin2} and \cite[Appendix, Sec.~3.3]{Tempel}).

An $\alpha$-invariant probability measure is called {\em ergodic\/} 
if every (measurable) invariant subset
of $\varOmega$ has either measure zero or measure one.  The ergodic
measures are exactly the extremal points of the convex set $\cPGO$. The
dynamical system $\Oomega$ is called {\em uniquely ergodic\/} if $\cPGO$
is a singleton set, i.e., if it consists of exactly one element.
As usual, $\Oomega$ is called {\em minimal\/} if, for all 
$\omega\in\varOmega$, the $G$-orbit
$\{\alpha^{}_t\ts \omega : t \in G\}$ is dense in $\varOmega$.
If $\Oomega$ is both uniquely ergodic and minimal, it is called
{\em strictly ergodic}.

\smallskip

Given an $m\in \cPGO$, we can form the Hilbert space $\LO$ of square
integrable measurable functions on $\varOmega$. This space is equipped
with the inner product
\begin{equation*}
    \langle f, g\rangle \; = \;
    \langle f, g\rangle^{}_\varOmega \; := \; 
    \int_\varOmega \overline{f(\omega)}\, g(\omega) \dd m(\omega).
\end{equation*}
The action $\alpha$ gives rise to a unitary representation $T =
T^\varOmega := T^{(\varOmega,\alpha,m)}$ of $G$ on $\LO$ by
\begin{equation*}
  T_t \! : \; \LO \; \longrightarrow \; \LO \, , 
  \quad (T_t f) (\omega) \; := \;
  f(\alpha^{}_{-t}\ts \omega) \, ,
\end{equation*}
for every $f\in \LO$ and arbitrary $t\in G$. 

An $f\in \LO$ is called an {\em eigenfunction\/} of $T$ with
{\em eigenvalue\/} $\hat{s}\in \Ghat$ if $T_t f = (\hat{s}, t) f$ for every
$t\in G$.  An eigenfunction (to $\hat{s}$, say) is called {\em continuous\/} 
if it has a continuous representative $f$  with
$f(\alpha^{}_{-t} \ts \omega) = (\hat{s},t)\ts f(\omega)$, for all
$\omega\in\varOmega$ and $t\in G$. The representation $T$ is said to have 
{\em pure point spectrum\/} if the set of its eigenfunctions is total in $\LO$.  
One then also says that the dynamical system $\Oomega$ has
{\em pure point dynamical spectrum}.

By Stone's theorem, compare \cite[Sec.~36D]{Loomis}, there exists a projection
valued measure
\begin{equation*} 
  E_T\! : \; \mbox{Borel sets of $\Ghat$} \; \longrightarrow \;
  \mbox{projections on $\LO$} 
\end{equation*} 
with
\begin{equation*} 
  \langle f, T_t f \rangle \; = \; 
  \int_{\Ghat} (\hat{s},t) \dd 
  \langle f,E_T(.) f\rangle (\hat{s}) \; := \; 
  \int_{\Ghat} (\hat{s}, t) \dd \rho^{}_f (\hat{s})\, , 
\end{equation*} 
where $\rho^{}_f = \rho_f^\varOmega:=
\rho_f^{(\varOmega,\alpha,m)}$ is the measure on $\Ghat$ defined by
$\rho^{}_f (B) := \langle f, E_T (B)f\rangle$.  In fact, by Bochner's
theorem \cite{Rudin}, $\rho^{}_f$ is the unique measure on $\Ghat$ with $\langle f,
T_t f \rangle = \int_{\Ghat}\, (\hat{s}, t) \dd \rho^{}_f 
(\hat{s})$ for every $t\in G$.

\section{Factors} \label{factor}

Factors of dynamical systems and the corresponding subrepresentations will be
an important tool in our study of deformations.  In this section, we recall their
basic theory, most of which is  well known. Since
details are somewhat scattered in the literature, we sketch some of
the proofs for the sake of completeness, or give precise references. 
Readers who are familiar with it, or are more interested to first
learn about diffraction, may skip this section at first reading.

\smallskip

Let $\Oomega$ 
and $\Ttheta$ be two topological dynamical systems under the action
of $G$, with a mapping 
$\varPhi \! : \varOmega\longrightarrow \varTheta$ that
gives rise to the following diagram:
\begin{equation} \label{diagram1}
\begin{CD}
  \varOmega @>\alpha>> \varOmega  \\
   @V{\varPhi}VV       @VV{\varPhi}V  @. \\
   \varTheta @>\beta>>  \varTheta
\end{CD}
\end{equation}

\smallskip

\begin{definition}  
  Let two topological dynamical systems\/ $\Oomega$ and\/ $\Ttheta$
  under the action of $G$ and a  mapping\/ 
  $\varPhi \! :  \varOmega \longrightarrow\varTheta$ be given.
  Then,  $\Ttheta$ is called a\/ {\em factor} of $\Oomega$, with factor map\/ 
  $\varPhi$, if\/ $\varPhi$ is a continuous surjection that makes the
  diagram\/ $(\ref{diagram1})$ commutative, i.e., 
  $\varPhi (\alpha^{}_t (\omega)) = \beta^{}_t (\varPhi (\omega))$ 
  for all\/ $\omega\in \varOmega$ and\/ $t\in G$.  
\end{definition}

Factors inherit many features from the underlying dynamical
system. Due to the commutativity of diagram \eqref{diagram1}, a
period $t\in G$ of $\omega$ is also a period of $\varPhi(\omega)$.
Clearly, the converse need not be true, as we shall see in an
example later on. Let us next recall three other properties
of dynamical systems which are inherited by factors.

\begin{fact}\label{top} 
  Let\/ $\Ttheta$ be a factor of\/ $\Oomega$, with factor map\/
  $\varPhi \! : \varOmega\longrightarrow \varTheta$. 
  Then, $U\subset \varTheta$ is open if and only if\/
  $\varPhi^{-1} (U)$ is open in\/ $\varOmega$. 
\end{fact}
\begin{proof} 
As $\varPhi$ is continuous, the only if part is clear. 
So, assume that $U\subset \varTheta$ is given with $\varPhi^{-1} (U)$ open. 
Then, $\varPhi^{-1} (\varTheta \setminus U) = \varOmega \setminus 
\varPhi^{-1} (U)$ is closed and thus compact, as $\varOmega$ is compact. 
Thus, by continuity  and surjectivity of $\varPhi$, the set 
$\varTheta\setminus U= \varPhi (\varPhi^{-1} (\varTheta\setminus U))$ 
is compact and, in particular, closed. Thus, $U$ is open.  
\end{proof}

Clearly, $\varPhi$ induces a mapping $\varPhi_* \! : \, \cMO
\longrightarrow \cMT$, $\mu\mapsto\varPhi_*(\mu)$, via
$\big(\varPhi_*(\mu)\big) (g) := \mu(g\circ\varPhi)$ for all
$g\in C(\varTheta)$. If $\mu$ is a probability measure on $\varOmega$,
its image, $\varPhi_*(\mu)$, is a probability measure on $\varTheta$.
Moreover, if $\varPhi$ is a factor map, invariance under the group
action is preserved. So, in this case, we obtain the mapping
\begin{equation} \label{diagram2}
  \varPhi_* \! : \; \cPGO \; \longrightarrow \; \cPGT \, , \quad
  \mu \, \mapsto \, \varPhi_*(\mu) \, ,
\end{equation}
where we stick to the same symbol, $\varPhi_*$, for simplicity.

\begin{fact} \label{fact2}
  Let\/ $\Ttheta$ be a factor of\/ $\Oomega$, with factor map\/ 
  $\varPhi \! : \varOmega \longrightarrow \varTheta$. Then, 
  $\varPhi_*$ of\/ $(\ref{diagram2})$ is a continuous surjection. 
  Moreover, it satisfies\/ $\varPhi_{\ast}
  \big(\sum_i c_i \mu_i\big) = \sum_i c_i \, \varPhi_{\ast} (\mu_i)$, 
  whenever\/ $\sum_i c_i \mu_i$ is a finite convex combination of
  measures\/  $\mu_i\in \cPGO$. Finally,
  $\varPhi_{\ast}$ maps ergodic measures to ergodic measures,
  and thus extremal points of\/ $\cPGO$ to
  extremal points of\/ $\cPGT$.
\end{fact}
\begin{proof}
By \cite[Prop.~3.2]{DGS}, the mapping $\varPhi_{\ast}$ is continuous, 
and by \cite[Prop.~3.11]{DGS}, it is onto. Direct calculations show
$\varPhi_{\ast} \big(\sum_i c_i \mu_i\big) = 
\sum_i c_i\, \varPhi_{\ast} (\mu_i)$ for
every finite convex combination $\sum_i c_i \mu_i$ of
measures in $\cPGO$. 

Let $\mu\in\cPGO$ be ergodic, i.e., any $\alpha$-invariant measurable subset
$A$ of $\varOmega$ satisfies either $\mu(A)=0$ or $\mu(A)=1$. Consider
$\nu:=\varPhi_\ast (\mu) \in \cPGT$, and let $B$ be a $\beta$-invariant
measurable subset of $\varTheta$, i.e., $\beta_t (B)=B$ for all $t\in G$.
Clearly, one has $\nu(B) = \mu(\varPhi^{-1}(B))$, where $A:=\varPhi^{-1}(B)
=\{\omega\in\varOmega : \varPhi(\omega)\in B\}$ is $\alpha$-invariant, as
a consequence of \eqref{diagram1}. So, $\nu(B)=\mu(A)$ is either $0$ or $1$,
and $\nu$ is also ergodic. The final claim about the extremal points
is then standard, compare \cite[Prop.~5.6]{DGS}.
\end{proof}

\begin{fact} \label{transfer}
  Let\/ $\Ttheta$ be a factor of\/ $\Oomega$, with factor map\/ 
  $\varPhi \! : \varOmega\longrightarrow \varTheta$. 
  If\/ $\Oomega$ is uniquely ergodic, minimal or strictly ergodic, 
  the analogous property holds for\/ $\Ttheta$ as well. 
\end{fact}
\begin{proof} 
If $\Oomega$ is uniquely ergodic, $\cPGO$ is a singleton set,
and $\cPGT = \varPhi_*(\cPGO)$ must then also be a singleton
set, by Fact~\ref{fact2}. So, also $\Ttheta$ is uniquely ergodic.
Apparently, every $G$-orbit in $\varTheta$ is the image of 
a $G$-orbit in $\varOmega$, under the factor map $\varPhi$.  
Continuity of $\varPhi$ implies $\varPhi(C) \subset 
\varPhi\big(\overline{C}\big) \subset \overline{\varPhi(C)}$
for arbitrary $C\subset \varOmega$. If $C$ is dense, 
$\overline{C}=\varOmega$, and $\overline{\varPhi(C)}=\varPhi(\overline{C}) =
\varTheta$ because $\varPhi$ is onto. This shows that minimality is
properly inherited, and the last claim on strict ergodicity
is then obvious.
\end{proof}

Now, let $\Ttheta$ be a factor of $\Oomega$ with factor map $\varPhi \! : 
\varOmega\longrightarrow \varTheta$ and let $m\in \cPGO$ be
fixed. For the remainder of this section, we denote the induced
measure on $\varTheta$ by $n=\varPhi_\ast (m)$. Consider the mapping
\begin{equation}  \label{i-def}
   i^\varPhi \! : \; \LT \longrightarrow \LO \, , \quad
   f\mapsto f\circ \varPhi \, ,
\end{equation}
and let $p^{}_\varPhi \! :  \LO\longrightarrow \LT$
be the adjoint of $i^\varPhi$. The maps $i^\varPhi$ and 
$p^{}_\varPhi$ are partial isometries. More precisely, 
$i^\varPhi$ is even an isometric embedding because
\begin{equation*}
  \langle i^\varPhi (g), i^\varPhi (f) \rangle^{}_\varOmega
  \; = \; \int_{\varOmega} \overline{(g\circ\varPhi)}\, (f\circ\varPhi)
  \dd m \; = \;  \big(\varPhi_* (m)\big) ( \overline{g}f) \; = \;
  \langle g, f\rangle^{}_{\varTheta}
\end{equation*}
for arbitrary $f,g\in \LT$. As $i^\varPhi$ is an isometry from $ \LT$ 
with  range $i^\varPhi (\LT)$, standard theory of partial isometries 
(compare \cite[Thm.\ 4.34]{Wei})  implies 
\begin{equation*}
   p^{}_\varPhi \circ i^\varPhi \; = \; {\rm id}^{}_{\LT}
   \quad \mbox{and} \quad
   i^\varPhi \circ p^{}_\varPhi \; = \; P_{i^\varPhi (\LT)}\ts , 
\end{equation*}
where ${\rm id}^{}_{\LT}$ is the identity on $\LT$ and 
$P_{i^\varPhi (\LT)}$ is the orthogonal projection of
$\LO$ onto $\cV:=i^\varPhi (\LT)$. 
  
Given these maps, we can discuss the relation between the
spectral theory of $T^\varOmega$ and\/ $T^\varTheta$.  
\begin{theorem}  \label{thm1}
   Let\/ $\LO$ and\/ $\LT$ be the canonical Hilbert spaces attached
   to the dynamical systems\/ $\Oomega$ and\/ $\Ttheta$, with factor 
   map\/ $\varPhi$ and\/ $n=\varPhi_\ast (m)$. 
   Then, the partial isometries\/ $i_{}^{\varPhi}$
   and\/ $p^{}_{\varPhi}$ are compatible with the unitary 
   representations\/ $T^\varOmega$ and\/ $T^\varTheta$ 
   of\/ $G$ on\/ $\LO$ and\/ $\LT$, i.e.,
\begin{equation*} 
     i^\varPhi  \circ T_t^\varTheta \; = \;   
     T_t^\varOmega \circ i^\varPhi \quad \mbox{and}\quad 
     T_t^\varTheta \circ p^{}_\varPhi  \; = \; 
     p^{}_\varPhi \circ T_t^\varOmega \ts ,  
\end{equation*}
   for all\/ $t\in G$. 
   Similarly, the spectral families\/ $E_{T^\varTheta}$ and\/ 
   $E_{T^\varOmega}$ satisfy 
\begin{eqnarray*}
  i^\varPhi\circ E_{T^\varTheta} (\cdot)  \; = \;
   E_{T^\varOmega} (\cdot)  \circ  i^\varPhi 
   \quad \mbox{and}\quad 
   E_{T^\varTheta} (\cdot)\circ  p^{}_\varPhi      
   \; = \;  p^{}_\varPhi \circ E_{T^\varOmega} (\cdot) \ts . 
\end{eqnarray*}
  The corresponding measures satisfy\/ $\rho_g^\varTheta = 
  \rho_{i^\varPhi (g)}^\varOmega$ for every\/ $g\in \LT$.
\end{theorem}
\begin{proof}  
Let $g\in \LT$ be given. As $\varPhi$ is a factor map, 
a short calculation gives
\begin{equation*}
  \big(T_t^\varOmega(i_{}^{\varPhi} (g))\big) (\omega) \; = \;
  g( \varPhi (\alpha_{-t} \omega)) \; = \;
  g( \beta_{-t} \varPhi (\omega))  \; = \;
  \big((i_{}^{\varPhi} T_t^\varTheta) (g)\big) ( \omega)
\end{equation*}  
and the first of the equations stated above follows. 
The second follows by taking adjoints. 

Choose $g\in\LT$. As discussed above, $\rho_g^\varTheta$ is the unique 
measure on $\Ghat$ with
\begin{equation*}
  \langle g, T_t^{\varTheta} g \rangle^{}_\varTheta  \; = \; 
  \int_{\Ghat} (\hat{s}, t) \dd
  \rho_g^\varTheta (\hat{s})
  \, , \quad \mbox{for all }  t\in  G \ts . 
\end{equation*} 
Similarly, $\rho^{\varOmega}_{i^\varPhi (g)}$ is the unique 
measure on $\Ghat$ with
\begin{equation*}
  \langle i^\varPhi (g) , T_t^{\varOmega} i^\varPhi (g) 
  \rangle^{}_\varOmega \; = \;
  \int_{\Ghat} (\hat{s}, t) \dd 
  \rho_{i^\varPhi (g)}^\varOmega(\hat{s})
  \, , \quad \mbox{for all } t\in G \ts .
\end{equation*} 
Moreover, as $i^\varPhi$ is an isometry, we obtain from the 
statements proved so far that
\begin{equation*}
   \langle g, T_t^{\varTheta} g \rangle^{}_\varTheta \; = \;
   \langle i^\varPhi (g), i^\varPhi (T_t^{\varTheta} g) 
   \rangle^{}_{\varOmega}  \; = \;
   \langle i^\varPhi (g) , T_t^{\varOmega} i^\varPhi (g) 
   \rangle^{}_\varOmega \, .
\end{equation*}
Putting the last three equations together, we obtain
\begin{equation*}
\int_{\Ghat} (\hat{s}, t) \dd \rho_g^\varTheta (\hat{s}) = \int_{\Ghat} (\hat{s}, t) 
  \dd \rho_{i^\varPhi (g)}^\varOmega(\hat{s})
\end{equation*}
for every $t\in G$. By the mentioned uniqueness of the involved measures, 
this gives 
\begin{equation*} \rho_g^\varTheta = \rho_{i^\varPhi
  (g)}^\varOmega \ts .
\end{equation*}
This, in turn,  implies
\begin{equation*} 
 \langle g ,  E_{T^\varTheta} (B) g  \rangle^{}_{\varTheta}  
  \; = \; \rho^{\varTheta}_g (B) \; = \;  
  \rho^{\varOmega}_{i^\varPhi (g)} (B) \; = \; 
   \langle i^\varPhi (g), E_{T^\varOmega} (B)
   i^\varPhi (g)\rangle^{}_\varOmega \;  = \;
  \langle g,  p^{}_\varPhi
   E_{T^\varOmega} (B) i^\varPhi (g)\rangle^{}_\varOmega  
\end{equation*} 
for all Borel measurable $B\subset\Ghat$ and every $g\in \LT$. As $g\in
\LT$ is arbitrary,  we infer $
 E_{T^\varTheta} (\cdot)  =
  p^{}_\varPhi E_{T^\varOmega} (\cdot) i^\varPhi$.
\end{proof}

One succinct way to summarize the core of Theorem~\ref{thm1}
is to say that the following  diagram is commutative, with the
map $i_{}^{\varPhi}$ (resp.\ $p^{}_{\varPhi}$) being injective 
(resp.\ surjective).
\begin{equation} \label{diagram3}
\begin{CD}
  \LT @>{i_{}^{\varPhi}}>> \LO @>{p^{}_{\varPhi}}>>  \LT   \\
  @V{T^{\varTheta}}VV @V{T^{\varOmega}}VV @V{T^{\varTheta}}VV  \\
  \LT @>{i_{}^{\varPhi}}>> \LO @>{p^{}_{\varPhi}}>>   \LT 
\end{CD}
\end{equation}

\smallskip

\begin{coro}  
  Assume the situation of Theorem~$\ref{thm1}$ and
  define\/ $\cV= i_{}^{\varPhi} (\LT)$.
  Then, $U \! : \, \LT \longrightarrow \cV$, 
  $f\mapsto i_{}^{\varPhi} (f)$, is a unitary map, the subspace\/  
  $\cV$ of\/ $\LO$ is invariant 
  under\/ $T^\varOmega$, and the restriction\/ 
  $T^\varOmega|_{\cV} $ of\/ $T^\varOmega$ to\/ 
  $\cV$ is unitarily equivalent to\/ 
  $T^\varTheta$ via\/ $U$.
\end{coro}
\begin{proof} 
As $i_{}^{\varPhi}$ is an isometric embedding, the map 
$U \!:\ts \LT\longrightarrow i_{}^{\varPhi}(\LT)$ is unitary. 
By Theorem~\ref{thm1}, we have $i^\varPhi  \circ T_t^\varTheta 
= T_t^\varOmega \circ i^\varPhi$. Consequently, the space 
$\cV=i_{}^{\varPhi}(\LT)$ is invariant under
$T^\varOmega$, with $T^\varOmega|_{\cV} U g = U
T^\varTheta g$ for every $g\in \LT$.
\end{proof}

The foregoing results describe the relationship between
$T^\varTheta$ and $T^\varOmega$ in the general case. 
In the special case of pure point spectrum, we can be more explicit as
follows.
\begin{prop} \label{pointspectrumfactor}
  Let\/ $\Ttheta$ be a factor of the dynamical system\/ 
  $\Oomega$, with factor map\/ 
  $\varPhi \! : \varOmega\longrightarrow \varTheta$. 
  Let\/ $m\in \cPGO$ be given, $n=\varPhi_*(m)$, and let\/
  $\LT$ and\/ $\LO$ be the corresponding Hilbert spaces. 
  Then, the following assertions hold.
\begin{itemize}
\item[\rm (a)]
  If\/ $g$ is an eigenfunction of\/ $T^\varTheta$ to the eigenvalue\/ $\hat{s}$, 
  $i^\varPhi (g) = g\circ \varPhi$ is an eigenfunction of\/ $T^\varOmega$ to the 
  eigenvalue\/ $\hat{s}$. 
\item[\rm (b)]
  If\/ $T^\varOmega$ has pure point dynamical spectrum,
  the same is true of\/ $T^\varTheta$.
\end{itemize}
\end{prop}
\begin{proof} 
(a): Let $g$ be an eigenfunction of $T^\varTheta$. Then, 
$i^\varPhi (g) = g\circ \varPhi$ is an 
eigenfunction of $T^\varOmega$, as $i^\varPhi  \circ T_t^\varTheta 
=  T_t^\varOmega \circ i^\varPhi$ by Theorem~\ref{thm1}.\\ 
(b): If $T^\varOmega$ has pure point dynamical spectrum, there exists an
orthonormal basis of $\LO$ which entirely consists of eigenfunctions of
$T^\varOmega$. Now, by Theorem~\ref{thm1}, we have $T_t^\varTheta
p^{}_\varPhi  =  p^{}_\varPhi T_t^\varOmega$. Therefore, 
$p_\varPhi f$ is an eigenfunction of $T^\varTheta$ (or zero) if 
$f$ is an eigenfunctions of $T^\varOmega$. As $p^{}_\varPhi$ is 
onto, the statement follows. 
\end{proof}

Let us now discuss the {\em continuity}\/ of eigenfunctions. 
Recall that a sequence $(B_n)_{n\in \NN}$ of compact sets in $G$ with 
non-empty interior is called {\em van Hove}, if it exhausts $G$ and if
\begin{equation*}
   \lim_{n\to \infty} \frac{|\partial^{K} B_n|}{|B_n|} \; = \; 0
\end{equation*}
for every compact $K$ in $G$, where $\partial^K B :=\overline{((B+ K) \setminus B)} 
\cup ((\overline{G\setminus B} - K)\cap B)$.  

For $G= \RR^d$ and $G=\ZZ^d$, the following lemma (and much more) was shown by 
Robinson in \cite{Rob}. His proof carries over easily to our situation. 
For the convenience of the reader, we include a brief discussion.  

\begin{lemma}  \label{transfer-continuous}
   Let\/ $\Oomega$ be a uniquely ergodic dynamical system. Denote the unique 
   invariant probability measure on\/ $\varOmega$ by\/ $m$.  Let\/ $\hat{s}$ 
   be an eigenvalue of\/ $T=T^{(\varOmega,\alpha,m)}$. Then, the following 
   assertions are equivalent. 
\begin{itemize}
\item[(i)] There exists a continuous eigenfunction\/ $f$ to\/ $\hat{s}$
    $($i.e., $f$ is continuous  with\/ $f(\alpha^{}_{-t} (\omega)) = 
   (\hat{s}, t) f(\omega)$ for all\/ $t\in G$ 
    and\/ $\omega\in \varOmega )$. 
\item[(ii)] The sequence\/ $A_{B_n} (h)$ of continuous functions on\/ $\varOmega$, 
    defined by  
\begin{equation*} 
    A_{B_n} (h)(\omega) \; := \; \frac{1}{|B_n|} \int_{B_n} 
    \overline{(\hat{s},t)}\ts h(\alpha_{-t} (\omega)) \dd t \ts ,
\end{equation*}
    converges uniformly, for every van Hove sequence\/ $(B_n)$ and every\/
    $h\in C(\varOmega)$.  
\end{itemize}
\end{lemma} 
\begin{proof} 
(i) $\Longrightarrow$ (ii) (cf. \cite{Rob}).
If $f$ is the continuous eigenfunction, $\lvert f \rvert$ is  invariant  and
continuous. As $\Oomega$ is uniquely ergodic, we may assume, without 
loss of generality, that $|f(\omega)|=1$ for every 
$\omega \in \varOmega$. Let $h\in C(\varOmega)$ be given. Apparently, the function 
$g =  h \overline{f} $ is continuous. Therefore, by unique ergodicity, the functions
\begin{equation*}
   \frac{1}{|B_n|} \int_{B_n} g(\alpha_{-t} (\omega)) \dd t \; = \; 
   \frac{1}{|B_n|} \int_{B_n} h(\alpha_{-t} (\omega)) \ts
   \overline{f(\alpha_{-t}(\omega))} \dd t  \; = \;  
   \frac{\overline{f(\omega)}}{|B_n|} \int_{B_n} 
   \overline{(\hat{s},t)}\, h(\alpha_{-t} (\omega)) \dd t
\end{equation*}
converge uniformly in $\omega\in \varOmega$. Multiplying by $f$ and using 
$f\ts \overline{f} = 1$, we infer (ii). 

\smallskip \noindent
(ii) $\Longrightarrow$ (i). As $\hat{s}$ is an eigenvalue of $T$, the projection 
$E(\{\hat{s}\})$ onto the eigenspace of $\hat{s}$ is not zero. Since $C(\varOmega)$ 
is dense in $\LO$, there exists an $h\in C(\varOmega)$ with $E(\{\hat{s}\}) h \neq 0$. 
Now, by the von Neumann ergodic theorem, see \cite[Thm.\ 6.4.1]{Krengel} for a
formulation that allows its derivation in the generality we need it here,
the sequence $A_{B_n} (h)$ converges in $\LO$ to $E(\{\hat{s}\}) h$.  By assumption 
(ii), this sequence converges uniformly to a function $g$. Thus, $g= E(\{\hat{s}\}) h$ 
in $\LO$. Moreover, by uniform convergence, $g$ is continuous and satisfies 
$g(\alpha^{}_{-t} (\omega)) = (\hat{s},t) g(\omega)$ for every $\omega\in \varOmega$ and 
$t\in G$. This gives (i). 
\end{proof}

Lemma~\ref{transfer-continuous} has the following interesting consequence. 
\begin{prop} \label{CEF} 
   Let\/ $\Oomega$ be a uniquely ergodic dynamical system, all eigenfunctions 
   of which are continuous. If\/ $\Ttheta$ is a factor of\/ $\Oomega$, with factor 
   map\/ $\varPhi$, it is a uniquely ergodic dynamical system, all eigenfunctions 
   of which are continuous as well. 
\end{prop}
\begin{proof} Fact~\ref{transfer} gives that $\Ttheta$ is uniquely ergodic.
Let $\hat{s}$ be an eigenvalue of $T^\varTheta$. Then, $\hat{s}$ is an eigenvalue 
of $T^{\varOmega}$ by Proposition~\ref{pointspectrumfactor}. We now apply 
Lemma~\ref{transfer-continuous} to infer continuity.
To that end, choose an arbitrary $g\in C(\varTheta)$, wherefore
$h = g\circ \varPhi$ belongs to 
$C(\varOmega)$. By (i) $\Longrightarrow$ (ii) of Lemma~\ref{transfer-continuous}, 
the sequence $(A_{B_n} (h))$ converges uniformly for every van Hove sequence $(B_n)$. 
A short calculation then gives
\begin{equation*}
    A_{B_n} (h) (\omega) \; = \; 
   \frac{1}{|B_n|} \int_{B_n} \big(g\circ \varPhi\big) (\alpha_{-t} (\omega)) 
   \overline{(\hat{s},t)} \dd t \; = \; 
   \frac{1}{|B_n|} \int_{B_n} g(\beta_{-t} (\varPhi (\omega))) \ts
   \overline{(\hat{s},t)} \dd t.
\end{equation*}
As $\varPhi$ is onto, this shows uniform convergence of 
$\theta \mapsto \frac{1}{|B_n|} \int_{B_n} g(\beta_{-t} (\theta)) 
\overline{(\hat{s},t)} \dd t.$
As $g\in C(\varTheta)$ was arbitrary, this gives the desired continuity 
statement, by (ii) $\Longrightarrow$ (i) of Lemma~\ref{transfer-continuous}. 
\end{proof}

Although we have not made use of it so far, it is possible to express 
$p^{}_\varPhi$ via a disintegration. Since it is instructive and also
useful in applications, we finish this section by giving the details
for the case when $\varOmega$ and $\varTheta$ are metrizable. 
We are in the somewhat simpler situation that a continuous map
$\varPhi \! : \, \varOmega \longrightarrow \varTheta$ exists.
By standard theory, compare \cite[Thm.~5.8]{Furstenberg} and 
\cite[Thm.~4.5]{N}, there exists a measurable map
\begin{equation*} 
  k \! : \; \varTheta \; \longrightarrow \; \mathcal{M}(\varOmega)
  \, , \quad \vartheta \, \mapsto \, k^\vartheta 
\end{equation*} 
that satisfies the following three properties.
\begin{enumerate}
\item For $n$-almost every $\vartheta\in \varTheta$,
  $k^\vartheta$ is a probability measure on $\varOmega$  supported in
  $\varPhi^{-1} (\vartheta)$.
\item For all $f\in L^1(\varOmega,m)$, the function 
  $f^{\{k\}} \! :  \varTheta \longrightarrow \CC$, $\vartheta
  \mapsto k^\vartheta (f)$,  is integrable with respect
  to $n=\varPhi_\ast (m)$.
\item For all $f\in L^1(\varOmega,m)$, one has
  $n (f^{\{k\}} ) = m(f)$. 
\end{enumerate}
In terms of integrals, the last property reads
\begin{equation*}
   \int_{\varTheta} f^{\{k\}} \dd n \; = \; 
   \int_{\varTheta} \int_{\varPhi^{-1}(\vartheta)}
   f(\omega) \dd k^\vartheta (\omega) \dd n (\vartheta)
   \; = \; \int_{\varOmega} f \dd m \ts .
\end{equation*}

\smallskip \noindent
{\sc Remarks}. (1) Note that \cite[Thm.~4.5]{N} only deals with bounded functions $f$. 
However, using standard monotone class arguments, it is not hard to extend the statements
given there to functions $f\in L^1 (\varOmega,m)$. This yields (2) and (3).  \newline
(2) The function $f^{k}$ can also be considered as a  conditional expectation of $f$ 
(see part (i)  of \cite[Thm.~5.8]{Furstenberg}  or part (b) of  \cite[Thm.~4.5]{N}).

\medskip

Given this disintegration, one can now describe the action of 
$p^{}_\varPhi$ on $f\in\LO$ explicitly, namely in terms of 
partial averages over the fibres $\varPhi^{-1}(\vartheta)$. 
\begin{prop}
   Assume that $\varOmega$ and\/ $\varTheta$ are compact
   metric spaces, and let\/ $n=\varPhi_* (m)$ as before. 
   Then, the equation\/ 
   $\big(p^{}_\varPhi (f)\big) (\vartheta) = k^\vartheta (f)$ 
   holds for all\/ $f\in\LO$ and\/ $n$-almost every 
   $\vartheta \in \varTheta$.
\end{prop}
\begin{proof} 
Fix $f\in \LO$, and let $g\in\LT$ be arbitrary.  Then, $f$ belongs to 
$L^1 (\varOmega, m)$,  since $\varOmega$ is compact and  
$\overline{g\! \circ\! \varPhi}\cdot f$ belongs to $L^1 (\varOmega,m)$, 
as it is the product of two $L^2$ functions.  Using the properties of $k$,
 we can then calculate 
\begin{eqnarray*}
   \langle g, f^{\{k\}} \rangle^{}_{\varTheta} & = & 
   \int_{\varTheta}\, \overline{g} \ts f^{\{k\}} \dd n \;\, = \;\,
   \int_{\varTheta}\, \overline{g(\vartheta)}\,
   \int_{\varPhi^{-1}(\vartheta)} f(\omega) \dd k^\vartheta (\omega) 
   \dd n (\vartheta)  \\ 
   & = & \int_{\varTheta} \int_{\varPhi^{-1}(\vartheta)}
   \overline{g(\varPhi(\omega))} \, f(\omega)
   \dd k^\vartheta (\omega) \dd n (\vartheta)  \;\, = \;\, 
   n \Big( \big(\ts \overline{i^{\varPhi}(g)}\ts f \big)^{\{k\}} \Big) \\ 
   & = & m \big(\ts \overline{i^{\varPhi}(g)}\ts f \big)
   \;\, = \;\,  \int_{\varOmega}
   \overline{g(\varPhi(\omega))}\ts f(\omega) \dd m(\omega) 
   \;\, = \;\, \langle i^\varPhi(g), f \ts \rangle^{}_\varOmega \\[1mm]
   & = &  \langle g , p^{}_\varPhi(f)\rangle^{}_\varTheta \, .
\end{eqnarray*}
As $g\in \LT$ is arbitrary, this gives $f^{\{k\}} = p^{}_\varPhi (f)$ in
$\LT$, and our claim follows.
\end{proof}

\section{Diffraction theory of measure and Delone dynamical systems} 
\label{Measure}

In this section, we specify the dynamical systems we are dealing with
and discuss the necessary background from diffraction theory.  The
material is taken from \cite{BL}, where the proofs and further details
can be found. For related material dealing with point dynamical systems, 
we refer the reader to \cite{Dworkin,Hof,LMS-1,Martin2,Boris,Boris2}.

As discussed in the introduction, our main focus is on measure dynamical 
systems which includes the case of point dynamical systems. For the 
convenience of the reader, however, we start this section with a short 
discussion of point dynamical systems and discuss the 
general case of measures only afterwards. 

\medskip
Let $V$ be an open neighbourhood of $0$ in $G$. A subset $\varLambda$ of
$G$ is called {\em $V\!$-discrete\/} if every translate of $V$ contains at most
one point of $\varLambda$. Such sets are necessarily closed. A set is 
{\em uniformly discrete\/} if it is $V\!$-discrete for some open neighbourhood
$V$ of $0$. The set of $V\!$-discrete point sets in $G$ is abbreviated as
$\mathcal{D}_V (G)$, while
the set of all uniformly discrete subsets of $G$ is
denoted by $\cUD$. The set $\cUD$ (and actually even the set 
$\mathcal{C}(G)$ of all closed subsets of $G$) can be topologized by 
a uniformity as follows. For $K\subset G$ compact and $V$ an open 
neighborhood of $0$ in $G$, we set
\begin{equation*}
   \UKV \; := \; \{(P_1,P_2) \in \cUD \times \cUD : 
   P_1 \cap K \subset P_2 + V \mbox{ and } 
   P_2 \cap K \subset P_1 + V\} \ts . 
\end{equation*}
It is not hard to check that $\{\UKV : K\mbox{ compact, } V
\mbox{ open with $0\in V$}\}$ generates a uniformity (see \cite[Ch.~6]{Kel} 
for basics about uniformities), and hence, via the neighbourhoods
\begin{equation*}
   \UKV (P) \; := \; \{ Q : (Q,P)\in \UKV\}
   \; , \quad P\in \cUD \ts ,
\end{equation*}
a topology on $\cUD$. This topology is called the {\em local rubber
topology\/} (LRT). For each open neighbourhood $V$ of $0$ in $G$, the 
set $\mathcal{D}_V (G) $ is compact in LRT. Apparently, $G$
acts on $\cUD$ by translation. By slight abuse of notation,
this action is again called $\alpha$, i.e., we define
\begin{equation*}
  \alpha^{}_t (\varLambda) \; := \; \{ t + x : x \in \varLambda\} 
   \; = \; t + \varLambda \ts .
\end{equation*}
To distinguish (compact) sets of measures $\omega$ from sets of point 
sets $\varLambda$, we shall use the suggestive notation $\varOmega$ 
versus $\varOmega_{\rm p}$ from now on. 

\begin{definition} \label{pds}
  The pair\/ $\Oop$ is called a\/ {\em point dynamical system} 
  if $\varOmega_{\rm p}$ is a 
  closed\/ $\alpha$-invariant subset of\/ $\mathcal{D}_V (G)$ for a 
  suitable neighbourhood $V\!$ of\/ $0$ in $G$. 
\end{definition}

Apparently, every $\varLambda\in \cUD$ gives rise to a point dynamical system 
$(\varOmega(\varLambda),\alpha)$, where $\varOmega(\varLambda)$ is the closure of 
$\{\alpha_t (\varLambda) : t\in G\}$ in LRT and $\alpha$ is the  action induced 
from the natural action of $G$ on $\cUD$.

\smallskip
After this  short look at  point dynamical systems, we now introduce 
our main object of interest: measure dynamical systems. As mentioned already, 
they generalize point dynamical systems (see below for details). 

Let $C>0$ and a relatively compact open set $V$ in $G$ be given.  A
measure $\mu\in \MM $ is called {\em $(C,V)$-translation bounded\/} if 
$|\mu|(t + V) \leq C$ for all $t \in G$. 
It is called {\em translation bounded\/} if
there exists such a pair $C,V$ so that $\mu$ is $(C,V)$-translation 
bounded. The set of all $(C,V)$-translation bounded measures is 
denoted by $\MCV$, the set of all translation bounded measures by
$\MTB$. In the vague topology, the set $\MCV$ is a compact Hausdorff 
space. There is an obvious action of $G$ on $\MTB$, again denoted by $\alpha$,  given by
\begin{equation*}
   \alpha \! : \; G\times  \MTB \; \longrightarrow \; \MTB 
   \, , \quad (t,\mu) \, \mapsto \, \alpha^{}_t\ts \mu 
   \quad \mbox{with} \quad (\alpha^{}_t \ts \mu) \, := \,
   \delta^{}_t * \mu \ts .
\end{equation*}
Restricted to $\MCV$, this action is continuous.

Here, the convolution of two convolvable measures $\mu, \nu$ is
defined by 
\[
  (\mu*\nu) (\varphi) \; = \; \int_G \varphi(r + s) 
  \dd\mu(r) \dd\nu(s) \ts .
\]  

\begin{definition} 
  $\Oomega$ is called a dynamical system on the translation bounded 
  measures on\/ $G$\/ $(${\rm TMDS} for short\/$)$ if there exist a constant\/ $C >0$ and 
  a relatively compact open set\/ $V\subset G$ such that\/  
  $\varOmega$ is a closed $\alpha$-invariant subset of\/ $\MCV$.
\end{definition}

It is possible to consider a point dynamical system as a TDMS. 
Namely, define
\begin{equation*}
   \delta \! : \; \cUD \; \longrightarrow \; \MTB
   \, ,\quad \delta(\varLambda)  \, := \,
   \textstyle{\sum_{x\in \varLambda}}\, \delta_x \ts ,
\end{equation*} 
where $\delta_x$ is the unit point (or Dirac) measure at $x$. 
The mapping $\delta$ is continuous and injective.
\begin{lemma}\label{delta}  
   If\/ $\Oop$ is a point dynamical system, the mapping\/
   $\delta \! :\ts \varOmega_{\rm p}\longrightarrow 
   \varOmega := \delta(\varOmega_{\rm p})$ establishes 
   a topological conjugacy between the point dynamical system\/ 
   $\Oop$ and its image, the\/ {\rm TMDS}  $\Oomega$. 
\end{lemma}
\begin{proof}
By \cite[Lemma~2]{BL}, $\delta\!:\varOmega_{\rm p}
\longrightarrow\delta(\varOmega_{\rm p})$
is a homeomorphism that is compatible with the $G$-action
$\alpha$, i.e., $\delta\big(\alpha^{}_t (\varLambda)\big) =
\alpha^{}_t \big(\delta(\varLambda)\big)$ for all
$\varLambda\in\varOmega_{\rm p}$ and all $t\in G$. So, $\delta$
provides a topological conjugacy as claimed.
\end{proof}

Having introduced our models, we can now discuss some key issues
of diffraction theory. Let $\Oomega$ be a TMDS, equipped with an
$\alpha$-invariant measure $m\in \cPGO$. We shall need the mapping
\begin{equation*} 
   f \! : \; C_c (G) \longrightarrow C(\varOmega)
   \, , \quad f_\varphi (\omega) \, := \, 
   \int_G \varphi (-s) \dd\omega(s) \ts . 
\end{equation*}
Then, there exists a unique measure $\gamma=\gamma_m$ on $G$, called the
\textit{autocorrelation\/} (often called Patterson function in
crystallography \cite{Cowley}, though it is a measure in our setting) with
\begin{equation*}
   \gamma (\overline{\varphi} \ast \psi_{\!\_}) \; = \;
   \langle f_\varphi, f_\psi \rangle
\end{equation*} 
for all $\varphi,\psi \in C_c(G)$, where $\psi_{\!\_}(s) := \psi(-s)$.  
The convolution $\varphi \ast \psi$ is defined by
$(\varphi \ast \psi) (t) = \int \varphi (t- s ) \psi (s) \dd s$.
For a more explicit formulation in terms of a weighted average,
see \cite[Prop.\ 6]{BL}. 


The measure $\gamma$ is positive definite. 
Therefore, its Fourier transform is a positive
measure $\gammahat$; it is called the {\em diffraction measure}. 
This measure describes the outcome of a diffraction experiment,
see \cite{Cowley} for background material.

\smallskip

\noindent {\sc Remark}.
This concept of an autocorrelation is defined via the entire dynamical
system, which implicitly involves a local averaging procedure. The 
conventional approach uses a limit of a sequence of finite measures
along a van Hove averaging sequence in $G$. If the dynamical system is
(uniquely) ergodic, the two notions coincide \cite{BL}. In general, the
definition we use here has the advantage of removing the dependence of
the averaging sequence and automatically deals with the {\em typical}\/
autocorrelation, at least with reference to the measure $m$.

In view of the fact that, in reality, one always faces {\em finite}\/
structures, one can give a justification along the following lines. 
Among all elements of the full system that are compatible with a given 
finite part, ``typical'' ones are those to be considered, if no other
piece of information is available. This means to take into account
all structures which, after a small translation and/or up to
some tiny local deformation, coincide with a fixed finite patch.
One way to do so is to take an average over all these possibilities
(on the level of their autocorrelations), which is essentially what
our $\gamma$  does. In the situation
of unique ergodicity, compare \cite{BL}, the precise method for forming
the average is irrelevant -- the result is independent of it.

\begin{theorem}[Theorem 7 in \cite{BL}] \label{characterization}  
  Let\/ $\Oomega$ be a\/ {\rm TMDS} with invariant measure\/ $m$. 
  Then, the following assertions are equivalent.
\begin{enumerate}
\item[(i)]  The measure\/ $\gammahat$ is a pure point measure.
\item[(ii)] $T^\varOmega$ has pure point dynamical spectrum.  \qed
\end{enumerate}
\end{theorem}

\smallskip

Theorem~\ref{characterization} links pure point 
diffraction spectrum to pure point
dynamical spectrum. This is of particular relevance for our
considerations. It will allow us to set up a perturbation and
stability theory for pure point diffraction spectrum by studying
(perturbations of) dynamical systems. This is the abstract core of our
investigation, to be analyzed next.

\section{Deforming measure and Delone dynamical systems: Abstract setting}
\label{Framework}

In this section, we introduce a deformation procedure for dynamical
systems that, under certain conditions, is {\em isospectral}, i.e., 
the deformation does not change the dynamical spectrum. In particular,
we shall later consider deformations of regular model sets
and show that a relevant class of deformations preserves
their pure point diffraction property. As discussed
in the introduction, these considerations are motivated by questions
from the mathematical theory of quasicrystals. They generalize
the corresponding results in \cite{Hof3,BD}.

\smallskip

For pedagogic reasons, we start with a short discussion of deformations of 
Delone dynamical systems. This results in Theorem \ref{stabilitypoint}. The 
general case of measure dynamical systems is treated afterwards. 

\medskip

Let $\Oop$ be a Delone dynamical system, with $\varOmega_{\rm p}$ 
contained in $\cDV$, and consider a continuous mapping
$q \!:\ts \varOmega_{\rm p} \longrightarrow G$ whose image then is
a compact set. In fact, let us assume
that $q(\varOmega_{\rm p}) -q(\varOmega_{\rm p}) \subset V$
for some neighbourhood $V$ of $0\in G$. Note that there exists 
an open neighbourhood $V'$ of $0$ in $G$ with 
\begin{equation*}
 V' + q(\varOmega_{\rm p}) -q(\varOmega_{\rm p})  
 \; \subset \; V .
\end{equation*}
In particular, for arbitrary $\varLambda\in \varOmega_{\rm p}$ and $y,z\in
\varLambda$ with $y\neq z$, we have $y + q(\varLambda-y) \neq z + q(\varLambda-
z)$ as well as 
\begin{equation} \label{lam-q}
   \varLambda_q \; := \; 
   \{x + q(\varLambda -x) : x\in \varLambda \} 
   \; \subset \; \mathcal{D}_{V'} (G)\ts.
\end{equation}
$\varLambda_q$ can be viewed as a ``deformed'' version of $\varLambda$,
which exlains the terminology.
Moreover,  $\varOmega_{\rm p}^{q}:= \{\varLambda_q : \varLambda \in 
\varOmega_{\rm p}\}$ can rather directly be seen to be 
$\alpha$-invariant and closed in $\mathcal{D}_{V'}(G)$.
Thus,  $(\varOmega_{\rm p}^{q},\alpha)$ is a point dynamical system,
and we have a mapping $\varPhi^q\!:\, \varOmega^{}_{\rm p}\longrightarrow
\varOmega^q_{\rm p}$ given by $\varPhi^q(\varLambda) = \varLambda_q$. This 
map can easily be seen to be a factor map.

In fact, it turns out that we do not need $q$ to be defined on the whole of 
$\varOmega_p$ to obtain a factor map.  It suffices to have it defined on a 
``transversal''. To be more precise here, we introduce the following subset of
$\varOmega_{\rm p}$,
\begin{equation} \label{xi-def}
   \varXi \; := \; \{ \varLambda \in \varOmega_{\rm p} : 0 \in \varLambda\}\ts .
\end{equation}
Since the elements of $\varOmega_{\rm p}$ are non-empty point sets of $G$, 
it is clear that each $G$-orbit in $\varOmega_{\rm p}$ contains at least one
element of $\varXi$. Moreover, the following holds. 
\begin{lemma} \label{xi-is-compact}
  If\/ $\varOmega_{\rm p}$ is a point dynamical system under the action of the\/
  {\rm LCA} group\/ $G$, the subset\/ $\varXi$ of\/ \eqref{xi-def} is compact.
\end{lemma}
\begin{proof}
By definition,  $\varOmega_{\rm p}$ is a closed subset of $\cDV$ for a suitable 
neighbourhood $V$ of $0$ in $G$. As $\cDV$ is compact in LRT, $\varOmega_{\rm p}$  
is compact in LRT as well.  So, we need to show that $\varXi\subset\varOmega_{\rm p}$ 
is a closed set.

Let $(\varGamma_\iota)$ be a net in $\varXi$ (so, $0\in\varGamma_\iota$ for all 
$\iota$) which converges to some $\varLambda$, where the latter must then lie in
$\varOmega_{\rm p}$. Assume that $0\not\in\varLambda$. Since $\varLambda$ is
itself a closed subset of $G$, we know that $G\setminus\varLambda$ is an open
set. By assumption, this open set would contain $0$, and hence also an
entire open neighbourhood of $0$. This, however, contradicts the
convergence $\varGamma_\iota\longrightarrow\varLambda$ in the LRT.
\end{proof}

As $\varXi$ is compact, every continuous function $q$ on $\varXi$ can be extended 
to a continuous function $\tilde{q}$ on $\varOmega_{\rm p}$ (the latter being compact
and hence normal) by Tietze's extension theorem, compare \cite[Prop.\ 1.5.8]{Ped}. The 
very definition of $\varOmega_p^{\tilde{q}}$, compare \eqref{lam-q}, shows that it
only depends on $q$ (and not on the extension chosen). In this situation, we can 
thus consistently define
\begin{equation} \label{fun-extend}
   \varOmega_p^q  \;  := \; \varOmega_p^{\tilde{q}} \ts .
\end{equation} 

Now, we can state our result on Delone dynamical systems. 
\begin{theorem} \label{stabilitypoint}
  Let\/ $\Oop$ be a point dynamical system under the action of the
  LCA group\/ $G$ with $\varOmega_{\rm p}\subset \cDV$ for a
  suitable neighbourhood\/ $V$ of\/ $0$ in $G$.
  Let\/ $q\! :\, \varXi\longrightarrow G$ be
  continuous with\/ $q(\varOmega_{\rm p}) -q(\varOmega_{\rm p}) \subset V$.
  Then, the following assertions hold.
\begin{itemize}
\item[\rm (a)] If\/ $\Oop$ has pure point diffraction spectrum $($w.r.t.\ an
    invariant probability measure $m)$, 
    so does\/ $(\varOmega^{q}_{\rm p},\alpha)$ $($w.r.t.\ the measure\/ 
    $\varPhi^q_* (m))$.  
\item[\rm (b)]  If\/ $\Oop$ is minimal or uniquely ergodic, 
    then so is\/ $(\varOmega^{q}_{\rm p},\alpha)$.
\item[\rm (c)] If\/ $\Oop$  is uniquely ergodic with pure point 
    diffraction spectrum and all of its eigenfunctions are continuous, 
    the same holds for\/ $(\varOmega^{q}_{\rm p},\alpha)$. 
\end{itemize}
\end{theorem}

\begin{proof}  
Let $\tilde{q}$ be a continuous extension of $q$ from $\varXi$ to 
$\varOmega_{\rm p}$. As discussed above in \eqref{fun-extend}, we then 
have a factor map $\varPhi^q : \varOmega_{\rm p} \longrightarrow  
\varOmega^{q}_{\rm p}$.  Now, we can prove the assertions.

\smallskip

(a): If $\Oop$ has pure point diffraction spectrum, it has pure point
dynamical spectrum, by Theorem~\ref{characterization}.  As
$ (\varOmega^{q}_{\rm p},\alpha) $ is a factor of $\Oomega$, it has pure 
point dynamical spectrum as 
well, by Proposition~\ref{pointspectrumfactor}. Now, another application of 
Theorem~\ref{characterization} shows that $(\varOmega^{q}_{\rm p},\alpha)$
has pure point diffraction spectrum.

\smallskip

(b): This follows from Fact \ref{transfer}. 

\smallskip

(c): The  statement about continuity of the eigenfunctions is immediate from
Proposition~\ref{CEF}. The other statements follow by (a) and (b). 
\end{proof}

\medskip

Having discussed the special case of point dynamical systems, we now treat 
the general case.  Let $\Oomega$ be a TMDS. We shall deform $\Oomega$ by means 
of a measure-valued mapping
\begin{equation*}
  \lambda  : \; \varOmega \; \longrightarrow \; \MM \, , \quad
  \omega \, \mapsto \, \lambda^\omega ,
\end{equation*}
which satisfies the following two properties.
\begin{itemize}
\item[(D1)] The mapping $\varOmega\times C_c(G)\longrightarrow\CC$,
            $(\omega,\varphi)\mapsto \lambda^{\omega}(\varphi)$, is continuous. 
\item[(D2)] There exists a compact $K\subset G$ such that
  $\supp(\lambda^\omega)\subset K$ for all $\omega\in \varOmega$. 
\end{itemize}
Such a deformation map $\lambda$ will be called {\em admissible}.
This definition entails the case that $\lambda^\omega\equiv\delta^{}_0$,
which we shall call the trivial deformation map.

\begin{prop} \label{beschraenkt} 
  Let\/ $\Oomega$ be a\/ {\rm TMDS} and let\/ 
  $\lambda \! : \varOmega \longrightarrow \MM$ 
  be an admissible deformation map.  Then, 
  $\omega\mapsto |\lambda^\omega| (1)$ is bounded.
\end{prop}
\begin{proof} 
Let $K$ be given according to (D2), and let $V\subset G$ be open and 
relatively compact. Since $\overline{K+V}$ is compact, one has
\begin{equation*}
   \lvert\lambda^\omega\rvert (1) \; = \; 
   \lvert\lambda^\omega\rvert (K+V) \; = \; 
   \sup\{\lvert\lambda^\omega (\varphi)\rvert : 
   \supp(\varphi) \subset \overline{K + V}, \, 
   \|\varphi\|_\infty \le 1\} \ts ,
\end{equation*}
where we used \cite[Prop.~1]{BL} in the last step. Due to compactness
of $\varOmega$, the statement now follows from (D1) and the uniform 
boundedness principle (see \cite[Thm.\ 2.2.9]{Ped}).  
\end{proof}

For $\omega\in \varOmega$ and $\varphi \in C_c (G)$, we define the 
actual {\em deformation} of $\omega$ into $\varPhi^\lambda (\omega)$ via
\begin{equation*}
   \big(\varPhi^\lambda (\omega)\big)(\varphi) \; := \; 
   \int_G \int_G \varphi (r + s) \dd
   \lambda^{\alpha^{}_{-r} (\omega)} (s) 
   \dd \omega (r) \ts ,
\end{equation*} 
where the double integral exists by (D1) and (D2). The constant
deformation map $\lambda^\omega\equiv\delta^{}_t$, with $t\in G$,
results in a translation, i.e., $\varPhi^\lambda (\omega)
=\delta^{}_t * \omega$ in this case, for all $\omega\in\varOmega$.
The trivial deformation map thus induces the identity. 
In general, the following is true.

\begin{prop} \label{prop4} 
  Let a TMDS\/ $\Oomega$ be given and let\/ $\lambda$ be an admissible
  deformation map. Then, the following assertions hold.
\begin{itemize}
\item[\rm (a)] For every\/ $\omega\in \varOmega$, the map\/
    $\varPhi^\lambda (\omega) \! : \, C_c (G) \longrightarrow \CC$, 
    $\varphi \mapsto \big(\varPhi^\lambda (\omega)\big)(\varphi)$, 
    is continuous, i.e., $\varPhi^\lambda (\omega)$ belongs to $\MM$. 
    Moreover, the map\/ $\varPhi^\lambda\! : \, \varOmega
    \longrightarrow \MM$, $\omega\mapsto \varPhi^\lambda (\omega)$, is
    continuous as well.
\item[\rm (b)] There exists a constant\/ $C >0 $ and an open 
    neighbourhood\/ $V\!$ of\/ $0$ in\/ $G$ such that\/ 
    $\varPhi^\lambda (\omega)$ belongs to\/ 
    $\MCV$, for all $\omega\in \varOmega$.
\item[\rm (c)] For all\/ $t\in G$ and\/ $\omega \in \varOmega$, one 
    has $\varPhi^\lambda (\alpha_t(\omega)) = 
    \alpha_t (\varPhi^\lambda (\omega))$. 
\end{itemize}
\end{prop}
\begin{proof} 
(a): Let $K$ be compact according to (D2). For fixed
$\omega \in \varOmega$ and $\varphi \in C_c (G)$ with support in
the compact set $L$, the function
\begin{equation*} 
  r \, \mapsto \, \int_G \varphi (r + s) \dd
    \lambda^{\alpha^{}_{-r}(\omega)} (s)
\end{equation*}
has support contained in $L-K$. Moreover, this function is continuous,
since it can easily be expressed as a composition of continuous
functions. In fact, extending this type of reasoning, one can show
that
\begin{equation*}
   F \! : \; \varOmega \times C_c (G) \; \longrightarrow \; C_c (G)
   \, ,  \quad F(\omega,\varphi) (r) \, := \, \int_G \varphi (r + s)  
   \dd\lambda^{\alpha^{}_{-r}(\omega)} (s) \ts ,
\end{equation*}
is continuous. In particular, $C_c(G)\longrightarrow \CC$, $\varphi
\mapsto \omega (F(\omega,\varphi))$, is continuous for fixed $\omega\in \varOmega$
and $\varOmega\longrightarrow \CC$, $\omega\mapsto
\omega(F(\omega,\varphi))$, is continuous for $\varphi \in C_c
(G)$. As
$$ \varPhi^\lambda (\omega)(\varphi)= \omega(F(\omega,\varphi)),$$
we infer (a).\\
(b): Let $L$ be an arbitrary non-empty open
set with compact closure. Let $1_{\overline{L} - K}$ be the
characteristic function of $\overline{L} - K$, where $K$ is taken from
(D1).  Then, for every $\varphi \in C_L (G)$, we have
\begin{equation*} 
  |\varPhi^\lambda(\omega)(\varphi)| \; \le \; 
  \int_G \, |\lambda|^{\alpha^{}_{-r}\omega} (1) \,
  \|\varphi\|^{}_\infty \, 1^{}_{\overline{L} -K} (r) 
  \dd\lvert\omega\rvert(r) \; \le \; C(\lambda) 
  \|\varphi\|_\infty \, \lvert\omega\rvert \, (\ts \overline{L} - K)\ts ,
\end{equation*}
where $C(\lambda)$ is the bound on $\omega \mapsto
|\lambda^\omega|(1)$ obtained in Proposition \ref{beschraenkt}.  Thus,
\begin{eqnarray*}
   \lvert\varPhi^\lambda (\omega)\rvert\ts (L + t) & = & 
   \sup\ts \{ \lvert \varPhi^\lambda (\omega) (\varphi)\rvert  :  
   \varphi \in C^{}_{L+t} (G), \|\varphi\|_\infty \le 1 \} \\
   & \le & C(\lambda) \| \varphi \|_\infty \,
           \lvert\omega\rvert\, (t+\overline{L}-K)
\end{eqnarray*}
is uniformly bounded in $t\in G$, as $\omega$ is translation bounded,
and (b) follows. \\ 
(c): This is immediate from
\begin{eqnarray*}
   \big(\varPhi^\lambda (\alpha^{}_t\ts \omega)\big) (\varphi) &=& 
   \int_G \int_G \varphi (r + s) \dd \lambda^{\alpha^{}_{-r+t}(\omega)} 
   (s) \dd (\alpha^{}_t\ts\omega)(r) \\ &=&\int_G \int_G \varphi (r + s + t) 
   \dd\lambda^{\alpha^{}_{-r} (\omega)} (s) \dd\omega(r)\\ 
   &=&
\big(\alpha^{}_t (\varPhi^\lambda (\omega)) \big)(\varphi)\ts ,
\end{eqnarray*} 
which is valid for every $\varphi \in C_c (G)$. 
\end{proof}

Define the set of periods of a measure $\omega$ as
\begin{equation} 
    {\rm Per}(\omega) \; := \;
   \{ t\in G : \alpha^{}_t \ts\omega = \omega\}\ts .
\end{equation}
We then have the following consequence.
\begin{coro} \label{transfer-of-periods}
   Let\/ $\Oomega$ be given and let\/ $\lambda$ be an admissible
   deformation map. For any\/ $\omega\in\varOmega$, with resulting
   deformation\/ $\varPhi^{\lambda}(\omega)$, one has
\[
   {\rm Per}(\omega) \; \subset \; 
   {\rm Per}(\varPhi^{\lambda}(\omega)) \ts .
\]
   Moreover, if any\/ $\omega\in\varOmega$ exists where\/
   ${\rm Per}(\varPhi^{\lambda}(\omega))$ is a true superset of\/
   ${\rm Per}(\omega)$, the mapping\/ $\varPhi^{\lambda}\!:\,
   \varOmega\longrightarrow\MM$ fails to be injective.
\end{coro}
\begin{proof}
The first claim follows at once from part (c) of Proposition~\ref{prop4}.
For the second claim, let $t$ be a period of $\varPhi^{\lambda}(\omega)$
that is {\em not} a period of $\omega$. Then, $\omega \neq \alpha^{}_{t}\ts
\omega$, but their images under $\varPhi^{\lambda}$ are equal.
\end{proof}

Part (a) of Proposition~\ref{prop4} implies that, for a given TMDS 
$\Oomega$, the set
\begin{equation*}
   \varOmega^{\lambda} \; := \;
   \{\varPhi^\lambda (\omega) : \omega \in \varOmega\}
\end{equation*}
is compact, as it is the image of a compact set under a continuous map. Furthermore, by
part (c) of the same proposition, $\varOmega^{\lambda}$ is invariant under
$\alpha$. In fact, by part (b) of Proposition \ref{prop4}, $\varOmega^{\lambda}$ is a
subset of $\MCV$ for suitable $C,V$.  Putting this together, we have
proved the following result.

\begin{lemma}\label{factorlemma} 
  Let\/ $\Oomega$ be a\/ {\rm TMDS} and let\/
  $\lambda \! : \varOmega \longrightarrow \MM$
  be an admissible deformation map. 
  Then, $(\varOmega^{\lambda}, \alpha)$ is a\/ {\rm TDMS}. 
  Moreover, $(\varOmega^{\lambda}, \alpha)$ is a factor of\/ $\Oomega$, 
  with factor map\/ 
  $\varPhi^\lambda \! :\, \varOmega\longrightarrow \varOmega^{\lambda}$.
  \qed
\end{lemma}

If the situation of Lemma~\ref{factorlemma} applies, we call
$(\varOmega^{\lambda}, \alpha)$ an {\em admissible\/} deformation
of $\Oomega$, with deformation map $\lambda$.
The main abstract result of this paper now reads as follows. 

\begin{theorem} \label{stabilityabstract} 
  Let\/ $\Oomega$ be a\/ {\rm TMDS} and let\/
  $\lambda \! : \,\varOmega \longrightarrow \MM$ be an
  admissible deformation map. Then, the following assertions 
  hold.
\begin{itemize}
\item[\rm (a)] If\/ $\Oomega$ has pure point diffraction spectrum $($w.r.t.\ some
  invariant probability measure $m)$, 
  so does\/  $(\varOmega^{\lambda},\alpha)$ $($w.r.t.\ the corresponding induced measure\/$)$. 
\item[\rm (b)] If\/ $\Oomega$ is minimal or uniquely ergodic, 
   then so is\/ $(\varOmega^{\lambda},\alpha)$.
\item[\rm (c)] If\/ $\Oomega$  is uniquely ergodic with pure point 
   diffraction spectrum and all of its eigenfunctions are continuous, 
   the same holds for\/ $(\varOmega^{\lambda},\alpha)$. 
\end{itemize}
\end{theorem}

\begin{proof}  
The proof is essentially the same as the proof of Theorem~\ref{stabilitypoint}.

(a):  If $\Oomega$ has pure point diffraction spectrum, it has pure point
dynamical spectrum, by Theorem~\ref{characterization}.  As
$(\varOmega^{\lambda},\alpha)$ is a factor of $\Oomega$ by
Lemma~\ref{factorlemma}, it has pure point dynamical spectrum as 
well, by Proposition~\ref{pointspectrumfactor}. Now, another application of 
Theorem~\ref{characterization} shows that $(\varOmega^{\lambda},\alpha)$
has pure point diffraction spectrum.

\smallskip

(b): This follows from Fact \ref{transfer}. 

\smallskip

(c): The  statement about continuity of the eigenfunctions is immediate from
Proposition~\ref{CEF}. The other statements follow by (a) and (b). 
\end{proof}

\noindent {\sc Remarks}. (1) Of course, the previous discussion of TMDS includes 
the case of Delone dynamical systems treated at the beginning of the section. To 
see this, one has to apply the mapping $\delta \! : \, \varOmega_p \longrightarrow \varOmega$ 
introduced in the previous section.\\
(2) The discussion of point dynamical systems given above requires a non-overlapping
condition under deformation, here written as $q(\varOmega_{\rm p}) -
q(\varOmega_{\rm p})  \subset V$ for a suitable open set $V$. In
the TMDS setting, such a restriction is not necessary, which shows
once more the greater flexibility of the approach via measures.

\section{Model sets and their deformation}\label{Modelset} 

Model sets probably form the most important class of examples of aperiodic
order. In their case, one starts with a periodic structure in a high-dimensional 
space and considers a partial ``image'' in a lower dimensional
space. This image will not be periodic any more but still preserve many
regularity features due to the periodicity of the underlying high dimensional
structure.  For a survey and further references, we refer the
reader to \cite{Moody-old,Moody}.

\smallskip

Let us start with a brief recapitulation of the setting of a cut and
project scheme and the definition of a model set. We need two locally
compact Abelian groups, $G$ and $H$, where $G$ is also assumed to be
$\sigma$-compact, see \cite{Martin2} for the reasons why this is
needed. As usual, neutral elements will be denoted by $0$ (or by
$0_G, 0_H$, if necessary). A {\em cut and project scheme\/} emerges 
out of the following collection of groups and mappings:
\begin{equation} \label{candp}
\begin{array}{cccccl}
    G & \xleftarrow{\,\;\;\pi\;\;\,} & G\times H & 
        \xrightarrow{\;\pi^{}_{\rm int}\;} & H & \\
   \cup & & \cup & & \cup & \hspace*{-2ex} \mbox{\small dense} \\
    L & \xleftarrow{\; 1-1 \;} & \tilde{L} & 
        \xrightarrow{\,\;\quad\;\,} & L^\star & \\
   {\scriptstyle \parallel} & & & & {\scriptstyle \parallel} \\
    L & & \hspace*{-38pt}
    \xrightarrow{\hspace*{47pt}\star\hspace*{47pt}} 
    \hspace*{-38pt}& & L^\star
\end{array}
\end{equation}
Here, $\tilde{L}$ is a {\em lattice\/} in $G\times H$, i.e., a cocompact
discrete subgroup. The canonical projection $\pi$ is one-to-one between
$\tilde{L}$ and $L$ (in other words, $\tilde{L}\cap \{0_G\}\times H =\{0\}$),
and the image $L^\star = \pi^{}_{\rm int}(\tilde{L})$
is dense in $H$, which is often called the internal space. In view of
these properties of the projections $\pi$ and $\pi^{}_{\rm int}$, one
usually defines the $\star$-map as $(.)^\star\!: L \longrightarrow H$
via $x^\star := \big( \pi^{}_{\rm int} \circ (\pi|^{}_{\tilde{L}})^{-1}\big) (x)$,
where $(\pi|^{}_{\tilde{L}})^{-1} (x) = \pi^{-1}(x)\cap\tilde{L}$, for all $x\in L$. 

A {\em model set\/} is now any translate of a set of the form
\begin{equation}
   \oplam(W) \; := \; \{ x\in L : x^\star \in W \}
\end{equation}
where the {\em window\/} $W$ is a relatively compact subset of $H$ with
non-empty interior. Without loss of generality, we may assume that the
stabilizer of the window,
\begin{equation} \label{window-shifts}
  H_W \; := \; \{ c\in H : c + W = W \}\ts ,
\end{equation}
is the trivial subgroup of $H$, i.e., $H_W = \{0\}$. If this were not
the case (which could happen in compact groups $H$ for instance), one could
factor by $H_W$ and reduce the cut and project scheme accordingly \cite{Martin2,BLM}.
Furthermore, we may assume that $\langle W-W \rangle$, the subgroup of $H$
that is algebraically generated by the subset $W-W$, is the entire group,
i.e., $\langle W-W \rangle = H$, again by reducing the cut and project
scheme to this situation, compare \cite{Martin1} for details. 

There are variations on the precise requirement to $W$
which depend on the fine properties of the model sets one is interested
in, compare \cite{Moody,Martin2}. In particular, a model set is called 
{\em regular\/} if $\partial W$ has Haar measure $0$ in $H$, and 
{\em generic\/} if, in addition, $\partial W \cap L^\star = \varnothing$. 

As discussed immediately after Definition \ref{pds},  every  model set $\varLambda$ 
gives rise to the dynamical system $(\varOmega(\varLambda), \alpha)$. It is
one of the central results of this area, compare \cite{Moody,Martin2} and
references given there, that (regular) model sets provide a very natural generalization
of the concept of a lattice.

\begin{theorem} \cite{Martin2}\label{RegModSet}
   Regular model sets are pure point diffractive. In fact, 
   $(\varOmega(\varLambda),\alpha)$ is uniquely ergodic with pure point 
   dynamical spectrum and continuous eigenfunctions.  \qed
\end{theorem}

For our purposes, it is sufficient to restrict our attention to regular 
model sets where $W$ is a compact subset of $H$ with
$\overline{W^\circ} = W$ (in particular, $W$ then has non-empty 
interior and, due to regularity, a boundary of Haar measure $0$). This is 
motivated by the fact that
diffraction cannot distinguish two model sets $\oplam(W)$
and $\oplam(W')$ if the symmetric difference $W\triangle W'$ of the
windows has Haar measure $0$ in $H$.

\smallskip
A regular model set with compact window $W$ can be deformed as follows
\cite{Hof,BD}. Let $\vartheta\!: H\longrightarrow G$ be a continuous
function with compact support, which, in view of the discussion around
\eqref{fun-extend}, we may assume to include $W$ if necessary.
If $\varLambda = \oplam(W)$, one defines
\begin{equation} \label{deformed-model-set}
  \varLambda_{\vartheta} \; := \;
  \{ x + \vartheta(x^\star) : x \in \varLambda \} \; = \;
  \{ x + \vartheta(x^\star) : x \in L \mbox{ and } x^\star \in W \}\ts .
\end{equation}
To make sure that $\varLambda_\vartheta$ is still a Delone set, one
usually requires that the compact set $K:= \vartheta(H) - \vartheta(H)$
satisfies $K\subset V$ where $V$ is an open neighbourhood of $0\in G$
so that $\varLambda\in\mathcal{D}_V (G)$.

Note that $\varLambda_\vartheta$ (if it is Delone) has a well defined
density, and one obtains
\begin{equation} \label{density}
   {\rm dens}(\varLambda_\vartheta) \; = \; {\rm dens}(\varLambda)\ts .
\end{equation}
In other words, an admissible deformation does not change the density.

Our aim is now to show that the continuous mapping $\vartheta$ induces a 
deformation map $q$ on $\varXi$.  To do so, we shall need the following lemma. 
It essentially says that the $\star$-map on $\varLambda$ can be extended 
to a unique continuous map on $\varXi$.

\begin{lemma} \label{new-help}
  Let\/ $\varLambda = \oplam(W)$, with\/ $W=\overline{W^\circ}$ compact,
  be a regular model set and assume that\/ $H_W=\{0\}$.  Then, 
  the set\/ $\{\varLambda - x : x\in\varLambda\}$ is
  dense in the compact set\/ $\varXi$ and there is precisely one 
  continuous mapping\/
  $\sigma\!: \, \varXi\longrightarrow W$ with\/ 
  $\sigma (\varLambda - x) = x^\star$ for every\/ $x\in \varLambda$. 
\end{lemma}

\begin{proof} 
First, let us show that $\{\varLambda - x : x\in\varLambda\}$ is dense
in $\varXi$, the latter being compact by Lemma~\ref{xi-is-compact}.

\smallskip

To this end, let $\varGamma\in\varXi$ be given and consider an arbitrary
neighbourhood $U_{K,V} (\varGamma)$ of $\varGamma$, where $K\subset G$ is 
compact and $V$ is an open neighbourhood of $0$ in $G$.  Replacing $K$ by 
$K\cup \{0\}$ if necessary, we can assume  $0\in K$ without loss of
generality. We have to provide an element of the form $\varLambda - p$ with 
$p\in \varLambda$ which belongs to $U_{K,V} (\varGamma)$. 

To do so, choose a compact neighbourhood $V'$ of $0\in G$ with
\begin{equation*}
 V' + V' \subset  V \quad \mbox{and} \quad V' = - V'.
\end{equation*}
As $\varXi$ is a subset of $ \varOmega_p (\varLambda)$, which is the orbit closure 
of $\{t + \varLambda : t\in G\}$, there exists a $t\in G$ with
\begin{equation*}
   t + \varLambda  \, \in \, U_{K + V', V'} (\varGamma).
\end{equation*}
As $0$ belongs to both $\varGamma$ and $K$, we infer that
\begin{equation*}
   0 \, \in \,  \varGamma \cap K \, \subset \, 
   \varGamma \cap (K + V') \, \subset \, (t  + \varLambda) + V'. 
\end{equation*}
Therefore, $0 = t + p  + v'$ with $p\in \varLambda$ and $v'\in V'$, or, put differently, 
$p = - t  -  v' \in \varLambda$. This gives
\begin{equation*}
   \varLambda - p \; = \; \varLambda + t + v' 
   \, \in \, U_{K + V', V'} (\varGamma) + v'
   \, \subset \, U_{K,V} (\varGamma),
\end{equation*}
where the last inclusion follows by our choice of $V'$. As discussed above, 
this proves the density statement.

\smallskip

It remains to show the existence and uniqueness of a continuous map 
$\sigma \!:\, \varXi\longrightarrow H$ with 
$\sigma (\varLambda - x) = x^\star$ for every $x\in \varLambda$,
where the uniqueness will be an immediate consequence of the continuity of 
$\sigma$ and the already established denseness of 
$\{\varLambda - x : x\in\varLambda\}$ in $\varXi$. 

Existence: By \cite[Lemma~4.1]{Martin2}, for every $\varGamma\in \varXi$,  
the set 
\begin{equation} \label{def-sigma}
  \sigma(\varGamma)\; =\;  \bigcap_{y\in\varGamma} (W-y^\star)
\end{equation}
is a singleton set in $H$ (note that
the sign change in our formulation does not affect this statement). 
In the sequel, we shall tacitly identify the singleton set $\sigma(\Gamma)$ 
with its unique element. Then, $ \sigma$ can be considered as a map on 
$\varXi$ with values in $H$. 

By \eqref{def-sigma}, $\overline{\varGamma^\star} \subset  W - \sigma(\varGamma)$. As
$0\in\overline{\varGamma^\star}$, we infer  $0=w-\sigma(\varGamma)$
for some $w\in W$, and hence $\sigma(\varGamma)\in W$.
If $\varGamma = \varLambda - x$ for some $x\in\varLambda$, then we claim
that $x^\star \in \sigma(\varLambda - x) =
\bigcap_{y\in\varLambda - x} (W-y^\star)$. This is so because 
$y\in\varLambda - x$ implies $y = \ell - x$ for some $\ell\in\varLambda$,
hence $W-y^\star = W - (\ell^\star - x^\star) = (W-\ell^\star) + x^\star$.
Clearly, $\ell^\star\in W$, so $0\in W-\ell^\star$, and this gives
$x^\star\in W - y^\star$. With $y\in\varLambda - x$ arbitrary, we obtain
$\sigma(\varLambda - x) = \{x^\star\}$, as $\sigma(\Gamma)$ is a singleton set.

Next, following \cite[Prop.~4.3]{Martin2}, we can show continuity of the
mapping $\sigma$. Let $\varGamma\in\varXi$, and let $V = V(\sigma(\varGamma))$
be an open neighbourhood of $\sigma(\varGamma)$ in $H$. Since
$\sigma(\varGamma) = \bigcap_{y\in\varGamma} (W-y^\star)$ is a singleton set,
one has
\[ 
    \bigcap_{y\in\varGamma} (W-y^\star)\setminus V \; = \;
    \varnothing\ts .
\]
As $V$ is open, each $(W-y^\star)\setminus V$ is closed, hence also compact.
So, there must be a {\em finite} set $F\subset\varGamma$ such that we already
have $ \bigcap_{y\in F} (W-y^\star)\setminus V = \varnothing$. This implies
that a compact set $K$ exists such that
$ \bigcap_{y\in \varGamma\cap K} (W-y^\star)\setminus V = \varnothing$, so
\[
   \bigcap_{y\in \varGamma\cap K} (W-y^\star) \; \subset \; V \ts .
\]
This inclusion means that $\varGamma'\cap K = \varGamma\cap K$, for any
$\varGamma'\in\varXi$, implies $\sigma(\varGamma')\subset V$.
By a standard argument, this can now be turned into the claimed 
continuity of $\sigma$.
\end{proof}

We can now show how $\vartheta$ induces a deformation $q$. 
\begin{prop} \label{extend}
  Let\/ $\varLambda = \oplam(W)$, with\/ $W=\overline{W^\circ}$ compact,
  be a regular model set and  assume that\/ $H_W=\{0\}$. Let 
  $\vartheta\!: W\longrightarrow G$  be  continuous.  Then,  
  there is precisely one continuous mapping\/
  $q\!: \varXi\longrightarrow G$ with\/
  $q(\varLambda - x) = \vartheta(x^\star)$ for all\/
  $x\in\varLambda$.
\end{prop}
\begin{proof}
This follows directly from Lemma~\ref{new-help}: 
Uniqueness follows because $\{\varLambda - x : x\in\varLambda\}$ is dense
in $\varXi$. Existence follows as we can simply define 
$q :=\vartheta \circ \sigma$ with the $\sigma$ of Lemma~\ref{new-help}. 
\end{proof}

\medskip

\noindent {\sc Remark}.
Let us point out that continuity of $\vartheta$ is not necessary to obtain 
continuity of $\vartheta \circ \sigma$. In fact, it is easy to construct 
examples where $\vartheta$ may even have countably many points of 
discontinuity (at points of $L^\star$, in fact). 

\smallskip
Based on Proposition~\ref{extend}, we can now prove our result on 
deformed model sets. 
\begin{theorem} \label{thm6}
   Let\/ $\varLambda$ be a regular model set and $\vartheta : H \longrightarrow G$ 
   a continuous map.  Let\/ 
   $\varLambda_{\vartheta}$ be defined according to $\eqref{deformed-model-set}$,
   with the restriction that it is still a Delone set.
   Then, $\varLambda_{\vartheta}$ is pure point diffractive.  In fact, 
   the dynamical system\/  $(\varOmega(\varLambda_\vartheta), \alpha)$ is 
   uniquely ergodic with pure point dynamical spectrum and continuous eigenfunctions. 
\end{theorem}
\begin{proof} 
Consider the map  $q \! : \varXi \longrightarrow G$ constructed in 
Proposition~\ref{extend}.  Plugging in the definitions, we easily find 
$\varLambda_q = \varLambda_\vartheta$. This, in turn, gives
\begin{equation*}
    (\varOmega(\varLambda))^q \; = \;
    \varOmega(\varLambda_q) \; = \; 
    \varOmega(\varLambda_\vartheta)  \ts .
\end{equation*}
Thus, it suffices to show that $\big((\varOmega(\varLambda))^q,\alpha\big)$ is uniquely 
ergodic with pure point dynamical spectrum and continuous eigenfunctions. This, 
however, is immediate from Theorem~\ref{stabilitypoint}.
\end{proof}

\smallskip \noindent
{\sc Remark}. 
Let us mention that the abstract result of Theorem~\ref{thm6} has a very concrete
extension in that it is possible to calculate the diffraction of $\varLambda_{\vartheta}$
explicitly. For the Euclidean setting, this is explained in \cite{Hof3,BD},
and we illustrate it below in a concrete example.

\section{Example: The silver mean chain}  \label{example}

Let us explain the various notions with a simple example in one dimension,
compare \cite[Sec.\ 8.1]{BMinv}.
To this end, consider the two letter substitution rule
\begin{equation} \label{silver-mean-def}
    \sigma : \begin{array}{rcl}
    a & \mapsto & aba \\ b & \mapsto & a
    \end{array}
\end{equation}
which allows the construction of a bi-infinite (and reflection symmetric)
fixed point as follows.
Starting from the (admissible) seed $w^{}_1 = a|a$, where $|$ denotes the reference 
point, and defining $w^{}_{n+1} = \sigma(w^{}_{n})$, one obtains the
iteration sequence 
\begin{equation*}
   a|a \;\stackrel{\sigma}{\longmapsto}\; aba|aba
   \;\stackrel{\sigma}{\longmapsto}\; abaaaba|abaaaba
   \;\stackrel{\sigma}{\longmapsto}\; ... 
   \;\xrightarrow{\;n\to\infty\;}\; w = \sigma(w)
\end{equation*}
where $w$ is a bi-infinite word in the alphabet $\{a,b\}$ and 
convergence is in the obvious product topology as generated from
the alphabet together with the discrete topology.

The corresponding substitution matrix reads
\begin{equation*}
   M_{\sigma} \; = \; \begin{pmatrix} 2 & 1 \\ 1 & 0 \end{pmatrix}
\end{equation*}
where $M_{k\ell}$ is the number of symbols of type $\ell$ in the
word $\sigma(k)$, for $k,\ell \in \{a,b\}$.
This matrix is primitive, with Perron-Frobenius eigenvalue
$s = 1+\sqrt{2}$, which happens to be a Pisot-Vijayaraghavan number.
It is often called the silver mean, due to its continued fraction
expansion ($s=[2;2,2,2,\ldots]$, in contrast to 
$[1;1,1,1,\ldots] = (1+\sqrt{5}\,)/2$ for the golden mean).
The corresponding eigenvectors (left and right) code the
frequencies of the letters $a$ and $b$ in $w$, and also the
information for a proper geometric representation of $w$ as a
point set in $\RR$, such that the substitution turns into a
geometric inflation rule. One convenient choice here is to represent
$a$ by an interval of length $1+\sqrt{2}$, and $b$ by one of length $1$.
Their frequencies are $\frac{1}{2}\sqrt{2}$ and 
$\frac{1}{2}\big(2-\sqrt{2}\ts\big)$, respectively.

This is an example of a so-called Pisot substitution with two symbols,
and the derived point set is known to be a regular model set
(with the projection scheme yet to be derived). Also, the fixed
point is a non-singular (or generic) member of the LI-class defined by it.
At the same time, it is a Sturmian sequence, and we could have started
with a concrete cut and project scheme (then with the compatibility with
the inflation to be established). We prefer the former possibility here,
as there is a rather elegant number theoretic formulation which we shall
now use.

Let $\varLambda_a$ and $\varLambda_b$ denote the left endpoints of the
intervals of type $a$ and $b$, with our reference point (formerly marked
by $|$) being mapped to $0$ in this process. Both point sets are subsets
of the $\ZZ$-module
\begin{equation*}
   \ZZ[\sqrt{2}\ts] \; := \; \{ m + n \sqrt{2} : m,n \in\ZZ\}
\end{equation*}
which happens to be the ring of integers in the quadratic field
$\QQ(\sqrt{2}\ts)$. There is one non-trivial algebraic conjugation
in this field, defined by ${}^\star\!:\sqrt{2}\mapsto -\sqrt{2}$, which 
maps $\ZZ[\sqrt{2}\ts]$ onto itself. This will take the r\^ole of the
$\star$-map in the cut and project scheme, which looks as follows.
\begin{equation*} 
\begin{array}{rcccccl}
   & \RR & \xleftarrow{\,\;\;\pi\;\;\,} & \RR \times \RR & 
        \xrightarrow{\;\pi^{}_{\rm int}\;} & \RR & \\
   \mbox{\small dense} \hspace*{-2ex}
   & \cup & & \cup & & \cup & \hspace*{-2ex} \mbox{\small dense} \\
   & \ZZ[\sqrt{2}\ts] & \xleftarrow{\; 1-1 \;} & \tilde{L} & 
        \xrightarrow{\; 1-1 \;} &\ZZ[\sqrt{2}\ts] & 
\end{array}
\end{equation*}
where $\tilde{L} = \{(x,x^\star) : x\in\ZZ[\sqrt{2}\ts] \}$ is a
(rectangular) lattice in $\RR^2$. In comparison to the standard
situation of model sets, compare \cite{Moody}, this cut and project
scheme is self-dual, see also \cite[p.~418]{Moody-old}. In particular,
the $\star$-map is then one-to-one on $\ZZ[\sqrt{2}\ts]$.

An explicit geometric realization of $\tilde{L}$ with basis vectors is
\begin{equation} \label{lattice}
   \tilde{L} \; = \; \Big\langle \binom{\sqrt{2}}{-\sqrt{2}} , 
   \binom{1}{1} \Big\rangle_{\ZZ}
\end{equation}
which has the nice property that we can directly work with the standard
Euclidean scalar product for our further analysis (rather than with the
quadratic form defined by the lattice).

In particular, we shall later also need the {\em dual lattice}
\begin{equation} \label{dual-lattice}
   \tilde{L}^* \; = \; \{y\in\RR^2 : xy\in\ZZ \mbox{ for all } 
   x\in\tilde{L} \} \; = \;
   \Big\langle \frac{1}{4}\binom{\sqrt{2}}{-\sqrt{2}} , 
   \frac{1}{2}\binom{1}{1} \Big\rangle_{\ZZ}
\end{equation}
(note the different star symbol), which has the projections
\[
   L^\circ \; = \; \pi(\tilde{L}^*) \; = \;
   \Big\{\frac{1}{2} \big( m + \frac{n}{\sqrt{2}} \big) : 
   m,n \in \ZZ \Big\} \; = \; \pi^{}_{\rm int}(\tilde{L}^*) \; = \;
   ( L^\circ)^\star \ts .
\]
Note that the $\star$-map is well defined on the rational span
of $L$ which includes $L^\circ$.

Let us continue with the construction of our model set.
By standard theory for the fixed point of a primitive substitution, the
sets $\varLambda_a$ and $\varLambda_b$ satisfy the equations 
\begin{eqnarray*}
   \varLambda_a & = & s \varLambda_a \;\overset{.}{\cup}\; 
               \big( s \varLambda_a + (1+s) \big) \;\overset{.}{\cup}\; 
               s \varLambda_b \\
   \varLambda_b & = & s \varLambda_a + s
\end{eqnarray*}
with $s=1+\sqrt{2}$ from above, and $\overset{.}{\cup}$ denoting the
disjoint union of sets. Under the $\star$-map followed by taking
the closure, one obtains a new set of equations for the windows
$W_a = \overline{\varLambda_a^\star}$ and $W_b = \overline{\varLambda_b^\star}$,
\begin{eqnarray*}
   W_a & = & s^{\star} W_a \;\cup\; \big( s^{\star} W_a + (1+s^{\star}) \big) 
             \;\cup\; s^{\star} W_b \\
   W_b & = & s^{\star} W_a + s^{\star}
\end{eqnarray*}
where $s^\star = 1-\sqrt{2}$ is less than $1$ in absolute value. This new
set of equations constitutes a coupled iterated functions system
that is a contraction. By standard Hutchinson theory, there is a unique pair
of compact sets $W_a$ and $W_b$ that solves this system, compare
\cite[Thm.~1.1 and Sec.~4]{BMinv} for details. It is easy to check that this 
solution is given by
\begin{equation} \label{windows}
  W_a \; = \; \big[\ts \tfrac{\sqrt{2}-2}{2},\tfrac{\sqrt{2}}{2} \ts\big] \; , \quad
  W_b \; = \; \big[-\tfrac{\sqrt{2}}{2},\tfrac{\sqrt{2}-2}{2} \ts\big]\ts .
\end{equation}
From here, one can also see that $W=W_a\cup W_b = 
\big[-\frac{\sqrt{2}}{2},\frac{\sqrt{2}}{2}\ts\big]$ is the window for the
full set $\varLambda = \varLambda_a \dot{\cup} \varLambda_b$, 
with $W = \overline{W^\circ}$.
Moreover, since $\pm 1/\sqrt{2}$ are not elements
of $\ZZ[\sqrt{2}\ts ]$, we see that $\varLambda = \oplam(W) =
\oplam(W^\circ)$, so that $\varLambda$ (and also $\varLambda_a$ and 
$\varLambda_b$) are regular, generic (or non-singular) model sets. The 
density of $\varLambda$ is ${\rm dens}(\varLambda) = 1/2$.

\smallskip

The deformation is now achieved by a suitable function 
$\vartheta\!: \RR\longrightarrow\RR$ which is continuous on $W$ and
vanishes on its complement. This is consistent with \eqref{fun-extend}
because the deformation rule \eqref{deformed-model-set} does not require 
the knowledge of $\vartheta$ for any value outside of $W$.
A simple but interesting candidate is
\begin{equation} \label{deform-fun-ex}
   \vartheta(y) \; = \; \begin{cases}
               \alpha \ts y + \beta, & y\in W \\
                 0 ,         & y\not\in W   \end{cases}
\end{equation}
with some constants $\alpha,\beta\in\RR$. For admissible values of
$\alpha$, the affine nature of $\vartheta$ on $W$ has
the effect of changing the relative length ratio of the
$a$ and $b$ intervals, with $\beta$ being a global translation. 
It is easy to check that the admissible values of $\alpha$ include
\[
    -1 \; < \; \alpha \; < \; 3 + \sqrt{2}
\]
which results in the ratio
\begin{equation} \label{ratio}
    \varrho \; = \; 
    \frac{{\rm length}(a^{}_{\vartheta})}
         {{\rm length}(b^{}_{\vartheta})}
    \; = \; 1 + \frac{1-\alpha}{1+\alpha}\sqrt{2} \ts .
\end{equation}
Here, we use $a^{}_{\vartheta}$ and $b^{}_{\vartheta}$ for
the intervals that result from the deformation \eqref{deform-fun-ex}.
For a given ratio, the parameter $\alpha$ is given by
$\alpha = (\sqrt{2} + 1 - \varrho)/(\sqrt{2} - 1 + \varrho)$.
We shall come back to this discussion in the next section.

\smallskip
Of particular interest is the fact that one does not only get the
theoretical result of pure point diffraction, but also an explicit
formula for the diffraction measure. A detailed account for its
calculation can be found in \cite{BD}, which can also be derived
explicitly via Weyl's lemma on uniform distibution, compare
\cite{Martin1,Moody2001} for a formulation of the latter in the 
context of model sets. The result is
\begin{equation} \label{diffraction-formula}
   \widehat{\gamma}^{}_{\varLambda_{\vartheta}} \; = \;
   \sum_{k\in L^\circ} \lvert A_{\vartheta}(k) \rvert^2 \, \delta_k
\end{equation}
where the so-called Fourier-Bohr coefficients (or diffraction
amplitudes) are given by
\begin{equation}  \label{amplitudes}
   A_{\vartheta}(k) \; = \;
   \frac{1}{2\sqrt{2}} \int_W
   e^{2\pi i (k^\star y - k\ts \vartheta(y))} \dd y
\end{equation}
for all $k\in L^\circ$, and $A_{\vartheta}(k) = 0$ otherwise. Note 
that $A_{\vartheta}(0) \equiv 1/2 = {\rm dens}(\varLambda)$ in agreement
with a previous remark.

To arrive at \eqref{diffraction-formula} and \eqref{amplitudes}, one
first shows that $A_{\vartheta}(k)$ must vanish for all $k\not\in L^\circ$,
which is part of \cite[Thm.\ 2.6]{BD}. Then, let $k\in L^\circ$, and 
consider the
points of $\varLambda$ in a (large) finite patch, e.g., in the interval
$B_r (0)$ of radius $r$ around $0$. We denote such a patch by
$\varLambda^{(r)}$ and set
\[
    \varLambda^{(r)}_{\vartheta} \; = \;
       \{ x + \vartheta(x^{\star}) : x \in \varLambda^{(r)} \}\ts .
\]
If we place unit point measures at the points of 
$\varLambda^{(r)}_{\vartheta}$, we obtain a finite measure whose
Fourier transform exists and reads
\[ 
   \sum_{x' \in \varLambda^{(r)}_{\vartheta}}
   e^{-2\pi i k x'} \; = \; 
   \sum_{x\in\varLambda_{}^{(r)}}
   e^{-2\pi i(kx + k \vartheta(x^\star))} \; = \;
   \sum_{x\in\varLambda_{}^{(r)}}
   e^{2\pi i(k^\star x^\star - k \vartheta(x^\star))}
\]
where the last step used the fact that
$e^{-2\pi i(kx + k^\star x^\star)}=1$ for $k\in L^\circ$ and
$x\in L$. Now, after dividing by the length of $B_r(0)$, one obtains
the coefficient $A_{\vartheta}(k)$ by taking the limit as
$r\to\infty$, which exists and gives \eqref{amplitudes} by
Weyl's lemma.

Let us also mention that, if we use the formulation via measures,  
the diffraction formula \eqref{diffraction-formula} remains valid for
{\em all} (continuous) functions $\vartheta$, not just for those
which preserve the Delone property.

For our special choice \eqref{deform-fun-ex}, one obtains
\begin{equation} \label{ampli2}
   A_{\alpha,\beta}(k) \; = \; e^{-2\pi i \beta k}\, \frac{\sin(z)}{2z}
   \Big|_{z = \pi (\alpha k - k^\star)\sqrt{2} }
\end{equation} 
for all $k\in L^\circ$.

\section{Topological conjugacy and further aspects}\label{Further}

In this section, we briefly comment on the question whether 
$(\varOmega^{\lambda},\alpha)$ is topologically conjugate to $\Oomega$.
A deformed model set need not be topologically conjugate to the
undeformed system. In our silver mean example, with the deformation function
$\vartheta$ of \eqref{deform-fun-ex}, we can find values of the scaling
parameter $\alpha$ where the factor becomes periodic, while $\varLambda$
itself (which corresponds to $\alpha=\beta=0$) is aperiodic. 
In such a case, in view of Corollary~\ref{transfer-of-periods}, we cannot have
topological conjugacy. Note that, in contrast to \cite{CS}, we do {\em not\/}
keep track of the type of the intervals here. If we did that
(e.g., by giving different weights to the points of $a$ and $b$ intervals), 
topological conjugacy would always be preserved under the deformation.

In particular, $\alpha = 1$ (which gives $\varrho = 1$)
results in $\varLambda_{\vartheta} = 2 \ZZ + \beta$.
Eq.~\eqref{diffraction-formula} then reduces to 
$\widehat{\gamma}^{}_{\varLambda_{\vartheta}} = \frac{1}{4} \delta^{}_{\ZZ/2}$,
as it has to. This is a concrete example of the phenomenon of an extinction 
rule, which can often be used to detect situations where topological
conjugacy fails. Here, by analyzing \eqref{ampli2} in detail, one finds that
the Fourier-Bohr spectrum 
\[
   \varSigma_{\alpha,\beta} \; := \; 
   \{ k \in \RR : A_{\alpha,\beta}(k) \neq 0 \}
\]
is independent of $\beta$, but depends on $\alpha$. Concretely, one has
\[
   \big\langle \varSigma_{\alpha,\beta} \big\rangle_{\ZZ}
       \; = \; \begin{cases}
          \frac{1}{2} \ZZ \ts , & \alpha = 1 \\
          L^\circ  \ts ,        & \text{otherwise}\ts .  \end{cases}
\]
Here, the $\ZZ$-span is needed because one can have systematic extinctions
also for $\alpha \neq 1$. This happens for $\alpha \in \QQ$ and for
$\alpha = 1 + r \sqrt{2}$ with $r\in\QQ$, through solutions of
$\sin (z)=0$ in \eqref{ampli2}. Such an extinction phenomenon is usually
linked to the existence of symmetries. In our case, for these special
values of $\alpha$, the point set $\varLambda_{\vartheta}$ admits an 
inflation symmetry, and the extinctions can be understood from 
that \cite{Benji}, see \cite{Benji2} for a general discussion.

Whenever $\alpha \neq 1$, the deformed model $\varLambda_{\vartheta}$ set is 
actually topologically conjugate to the original model set $\varLambda$, though
in general not via a local derivation rule, compare \cite{CS2} for a recent
clarification of the relation between these concepts.

\smallskip

Another interesting phenomenon is the appearance of periodic diffraction,
even if the underlying structure is non-periodic. For simplicity,
let us concentrate on the case $\beta=0$. Whenever $\varrho$ of \eqref{ratio}
is a rational number, $\varrho = p/q$ say with $p,q$ coprime, the set of
positions of $\varLambda_{\vartheta}$ is a subset of a lattice in $\RR$
(of period $\lambda = {\rm length}(a^{}_{\vartheta})/p =
{\rm length}(b^{}_{\vartheta})/q$). Consequently, by \cite[Thm.\ 1]{MB},
the diffraction measure of the corresponding Dirac comb is periodic,
with period $1/\lambda$. As the diffraction is also pure point, by our
Theorem~\ref{thm6}, it is of the form $\mu * \delta_{\ZZ/\lambda}$, where
$\mu$ is a finite positive pure point measure on $[0,1/\lambda)$.
Unless $\alpha = 1$, the Fourier-Bohr spectrum is dense in $\RR$, and
the underlying Dirac comb based on $\varLambda_{\vartheta}$ is not periodic.
So, in our example, failure of topological conjugacy coincides with the
existence of periods for $\varLambda_{\vartheta}$.

\smallskip

In the example, and also in our general discussion, we started from
a model set and constructed a deformation scheme. In general, a deformation
will not result in another model set, though its Fourier-Bohr
spectrum remains unchanged. The latter is of central importance for the
actual structure determination in crystallography, e.g., from a
diffraction experiment. It is often implicitly assumed that the
underlying structure is a model set, but our above analysis shows
that this need not be the case. An important
open question is thus how to effectively characterize model sets versus
deformed model sets by means of intrinsic properties, preferably by
easily accessible ones. Some first results can be infered from  \cite{BLM},
but more has to be done in this direction.

\bigskip
\subsection*{Acknowledgements}
It is our pleasure to thank Robert V.\ Moody and Lorenzo Sadun
for a number of very helpful discussions.
This work was supported by the German Research Council (DFG).

\end{document}